\documentclass[pdflatex,sn-basic,iicol]{sn-jnl}

\usepackage{graphicx}%
\usepackage{epstopdf}
\usepackage{multirow}%
\usepackage{amsmath,amssymb,amsfonts}%
\usepackage{amsthm}%
\usepackage{mathrsfs}%
\usepackage[title]{appendix}%
\usepackage{xcolor}%
\usepackage{textcomp}%
\usepackage{manyfoot}%
\usepackage{booktabs}%
\usepackage{algorithm}%
\usepackage{algorithmicx}%
\usepackage{algpseudocode}%
\usepackage{listings}%
\usepackage{natbib}
\usepackage{tikz}
\usepackage{float}

\usetikzlibrary{calc,positioning}

\theoremstyle{thmstyleone}%
\theoremstyle{thmstyletwo}%
\newtheorem{remark}{Remark}%

\theoremstyle{thmstylethree}%

\newtheorem{problem}{Problem}

\DeclareMathOperator\Span{span}

\begin{document}

\title{Active vibration control of nonlinear flexible structures via reduction on spectral submanifolds}


\author[1]{\fnm{Cong} \sur{Shen}}

\author*[1]{\fnm{Mingwu} \sur{Li}}\email{limw@sustech.edu.cn}

\affil*[1]{\orgdiv{Department of Mechanics and Aerospace Engineering}, \orgname{Southern University of Science and Technology}, \orgaddress{\city{Shenzhen}, \postcode{518055}, \country{China}}}


\abstract{Large amplitude vibrations can cause hazards and failure to engineering structures. Active control has been an effective strategy to suppress vibrations, but it faces great challenges in the real-time control of nonlinear flexible structures. Here, we present a control design framework using reductions on aperiodic spectral submanifolds (SSMs) to address the challenges. We formulate high-dimensional nonlinear optimal control problems to suppress the vibrations and then use the SSM-based reductions to transform the original optimal control problems into low-dimensional linear optimal control problems. We further establish extended linear quadratic regulators to solve the reduced optimal control problems, paving the road for real-time active control of nonlinear flexible structures. We demonstrate the effectiveness of our control design framework via a suite of examples with increasing complexity, including a finite element model of an aircraft wing with more than 130,000 degrees of freedom.}




\maketitle

\section{Introduction}\label{sec1}
Vibration is commonly observed in engineering applications and can cause hazards to engineering structures~\cite{waniCriticalReviewControl2022}. Appropriate vibration control designs are essential to suppress severe vibrations of the structures and hence improve their safety~\cite{davidwagg,xieStructuralControlVibration2020}. In particular, flexible structures are increasingly used in soft robotics~\cite{websterMechanicsPrecurvedTubeContinuum2009,russoContinuumRobotsOverview2023} and aerospace engineering~\cite{lee-glauserSatelliteActivePassive1996}. The increased flexibility can enhance the performance of robots or reduce the weight of aerospace structures but often leads to large amplitude vibrations. Given these vibrations are nonlinear and of high (or infinite) dimensionality, it is difficult to control them precisely and effectively~\cite{mohamedVibrationControlVery2005}. Therefore, vibration control of flexible structures has been an important engineering problem and attracted the interest of many researchers~\cite{ khouryRecentAdvancesVibration2013, ghaediInvitedReviewRecent2017,kumarReviewControllersStructural2023}.

Control strategies for suppressing vibration can be grouped as passive, active~\cite{liSelfPoweredActiveVibration2022}, semi-active~\cite{luGeneralMethodSemiActive2004}, and hybrid control ~\cite{xieStructuralControlVibration2020,khouryRecentAdvancesVibration2013}. Passive control integrates energy dissipation or transfer devices into the structures to tune the parameters of the structures, such as mass, damping, and stiffness, and hence suppress the vibration of main structures~\cite{shiPerformanceComparisonPassive2017,thenozhiAdvancesModelingVibration2013}. In contrast, active control feeds force input via actuators to control the response of the structures directly~\cite{lee-glauserSatelliteActivePassive1996,korkmazReviewActiveStructural2011}. In particular, it can adjust the control input based on the feedback of external excitation and the system's response. Semi-active control is the combination of passive and active control~\cite{khouryRecentAdvancesVibration2013} and hence can combine the robustness of the passive control and the high performance of the active control~\cite{zhangActivepassiveIntegratedVibration2017,dingVibrationControlFluidconveying2023}. In summary, active control can rapidly suppress the vibration and has been shown as the most commonly used and effective approach for mitigating the adverse impacts of dynamic loads on structures~\cite{waniCriticalReviewControl2022}. 

Several active vibration control algorithms have been proposed in the past few decades. Among them, there are linear control schemes such as Proportional-Integral-Derivative (PID) control~\cite{thenozhiStabilityAnalysisActive2014}, linear quadratic regulator (LQR) control \cite{shiPerformanceComparisonPassive2017} and its extension linear quadratic Gaussian (LQG) control \cite{hoActiveVibrationControl2007}. More advanced control algorithms such as model predictive control (MPC)~\cite{takamotoComprehensivePredictiveControl2022}, $\mathcal{H}_{\infty}$ control \cite{cancielloSelectiveModalControl2017}, and sliding mode control(SMC)~\cite{choiVibrationControlElectrorheological2007} were also applied for the purpose of vibration suppression. Recently, control designs via artificial intelligence, such as fuzzy control and neural networks, were proposed~\cite{kumarReviewControllersStructural2023, liControlMethodologiesVibration2018,kandasamyReviewVibrationControl2016}.

The above control algorithms face challenges when directly applied to high-dimensional nonlinear flexible structures. Indeed, the vibration of these structures changes dynamically because the external excitation can change rapidly~\cite{gaoVibrationAnalysisControl2021, SHLX201504001}. This requires the control algorithm to rapidly cope with the fast-changing vibration~\cite{kangundeReviewDronesControlled2021} and then achieve real-time control. Meanwhile, the equations of motion of flexible structures are in the form of ordinary differential equations (ODEs), obtained via numerical approximations such as finite element method (FEM)~\cite{jain2022compute}. To guarantee the accuracy of these approximations, the discrete model of flexible structures usually has many degrees of freedom, resulting in high-dimensional nonlinear systems. However, solving active vibration control for such high-dimensional nonlinear dynamical systems is time-consuming \cite{wanDatadrivenModelReduction2024}, constituting a bottleneck of the real-time control.

Model order reduction provides a powerful framework to address the above challenges. Specifically, a high-dimensional system is replaced by a low-dimensional reduced-order model (ROM), which enables the real-time control~\cite{besselinkComparisonModelReduction2013, sumanInvestigationImplementationModel2022}. For linear structural dynamics, the core ingredient of model reduction is to project the full state space onto some lower-dimensional subspace. Depending on the selection and construction of the subspace, modal truncation~\cite{xianminActiveVibrationController2002}, moment matching~\cite{vakilzadehVibrationControlMicroscale2020,gildinModelControllerReduction2009a}, Galerkin projection via proper orthogonal decomposition (POD) modes~\cite{banksProperOrthogonalDecompositionbased2002}, and balanced truncation \cite{mathewsOrderReductionClosedLoop2002,antoulasApproximationLargescaleDynamical2005,kingReducedOrderControllers2006} have been used to construct ROMs for linear structural dynamics.

The above linear reduction methods have also been applied to reduce nonlinear structural dynamics. For instance, balanced truncation has been applied to realize model reduction considering rigid-flexible coupling and nonlinear characteristics of high-speed spinning beams~\cite {zhouIntelligentVibrationControl2015}. However, they destroy the feature of invariant manifolds, and hence, their constructions are problem-dependent, and the resulting ROMs can still be high-dimensional~\cite{hallerNonlinearNormalModes2016a,touzeModelOrderReduction2021}.

Invariant-manifold-based nonlinear model reduction schemes, in contrast, preserve the invariance feature and can hence achieve exact model reduction~\cite{hallerNonlinearNormalModes2016a,touzeModelOrderReduction2021}, which is essential for flexible structures where we need to account for geometric nonlinearities~\cite{maSuspensionNonlinearAnalysis2024}. Nechak and Lyes \cite{nechakRobustNonlinearControl2022} proposed a center manifold-based approach to reduce the dimension of nonlinear systems and then design a sliding-mode controller to suppress the vibration induced by friction. A center manifold is not structurally stable because it is not persistent under the addition of damping. For structural systems with hyperbolic fixed points, nonlinear normal modes (NNMs) that are defined as invariant manifolds tangent to spectral subspaces are relevant for model reduction~\cite{shaw1999modal,touze2006nonlinear,Liu2019SimultaneousNF,touzeModelOrderReduction2021}. However, as shown in~\cite{hallerNonlinearNormalModes2016a}, infinitely many NNMs exist for a given spectral subspace. Consequently, reductions via NNMs are not crystal clear from a mathematical point of view.

Among these infinitely many NNMs, there is a unique, smoothest one under proper non-resonance conditions, as pointed out by Haller and Ponsioen \cite{hallerNonlinearNormalModes2016a}. They defined the unique invariant manifold as a \emph{spectral submanifold} (SSM), which provides a solid foundation for constructing invariant-manifold-based ROMs. Indeed, SSM-based model reduction has emerged as a powerful tool for model reduction of nonlinear systems, including structural~\cite{jain2022compute,liNonlinearAnalysisForced2022,liModelReductionConstrained2023}, fluid~\cite{cenedese2022data,kaszas2024capturing}, and fluid-structure interactions~\cite{liNonlinearModelReduction2023a,li2024data,Xu_Kaszás_Cenedese_Berti_Coletti_Haller_2024}. SSM-based ROMs enable efficient and even analytic extractions of backbone curves~\cite{breunung2018explicit}, forced responses under harmonic forcing~\cite{breunung2018explicit,li2024fast}, and bifurcations of periodic and quasi-periodic orbits~\cite{liNonlinearAnalysisForced2022}.

The previous studies on SSM-based model reduction mainly focus on free and forced vibration with periodic forcing. However, in the control setting, the input signal can be arbitrary. The SSM theory has recently been extended to general non-autonomous dynamical systems either weakly forced or slowly varying~\cite{Haller2024NonlinearMR}. This extension provides a rigorous framework for constructing SSM-based ROMs with control. Indeed, SSM-based ROMs with control have been constructed for the trajectory tracking of soft robotics in a data-driven setting~\cite{alora2023data,alora2023robust}.

Here, we aim to establish an SSM-based control design framework for the active control of high-dimensional nonlinear flexible structures in equation-driven settings, given FE models of structures are often available. Our control design framework differs from that of~\cite{alora2023data,alora2023robust} in two perspectives. Firstly, simulating the full system is not involved because we use the equations of motion of structures to perform model reduction directly. Secondly, our control design consists of two ROMs: one is nonlinear without control and hence can be simulated in an offline stage, and the other one is linear with control. As we will see, this separation enables us to apply well-established linear control theory, significantly simplifying control analysis and design tasks. In particular, we will transform the original high-dimensional nonlinear control problem into a low-dimensional linear control problem using our SSM-based model reduction. This is one of the key contributions of our study.

The rest of this paper is organized as follows. Sect.~\ref{sec:problem-formulation} presents the problem formulation of this study, where we first present a system setup, followed by a formulation of optimal control problem to suppress the large amplitude nonlinear vibration. In Sect.~\ref{sec:review-ssm}, we revisit the theory of SSM and present the SSM-based model reduction of unforced systems. Then, we extend the SSM to a time-dependent manifold. With SSM-based model reduction, we transform the original nonlinear control problem into a linear one in Sect.~\ref{sec:Controltransformation}. We further introduce two reduction schemes in Sect.~\ref{sec:ControlL} to reduce the high-dimensional linear control problem to a low-dimensional one to achieve real-time control. The reduced linear control problem is an extended linear quadratic optimal control problem, which can be solved via an extended linear-quadratic regulator control, as we detailed in Sect.~\ref{sec6:LQ}. In Sect.~\ref{sec:numerialExamples}, we demonstrate the effectiveness of our control design framework via a suite of examples with increasing complexity. We conclude this study in Sect.~\ref{sec:conclusion}.

\section{Problem formulation}
\label{sec:problem-formulation}
\subsection{System setup}
\label{sec:setup}
\begin{sloppypar}
We consider a nonlinear flexible structure subject to external forcing and control input. Without loss of generality, the equations of motion for the discrete model of the structure are in the form:
\begin{align}
\label{eq:eom-second-full}
&\boldsymbol{M}\ddot{\boldsymbol{x}}+\boldsymbol{C}_\mathrm{d}\dot{\boldsymbol{x}}+\boldsymbol{K}\boldsymbol{x}+\boldsymbol{f}(\boldsymbol{x},\dot{\boldsymbol{x}})=\epsilon \left(\boldsymbol{E}(t)+\boldsymbol{D}\boldsymbol{u}(t)\right),\,\, \nonumber\\
&\boldsymbol{x}(t_0)=x_0,\,\, \dot{\boldsymbol{x}}(t_0)=\dot{\boldsymbol{x}}_0,\,\, 0\leq\epsilon\ll1,
\end{align}
where $\boldsymbol{x}\in\mathbb{R}^n$ is a generalized displacement vector ($n\gg1$); $\boldsymbol{M}\in\mathbb{R}^{n\times n}$ is a positive definite mass matrix; $\boldsymbol{C}_\mathrm{d},\boldsymbol{K}\in\mathbb{R}^{n\times n}$ are damping and stiffness matrices; $\boldsymbol{f}(\boldsymbol{x},\dot{\boldsymbol{x}})$ is a $C^r$ smooth nonlinear function such that
$\boldsymbol{f}(\boldsymbol{x},\dot{\boldsymbol{x}})\sim \mathcal{O}(|\boldsymbol{x}|^2,|\boldsymbol{x}||\dot{\boldsymbol{x}}|,|\dot{\boldsymbol{x}}|^2)$; $\boldsymbol{E}$ represents a vector of external excitation to which the structure is subjected; $\boldsymbol{D}\in\mathbb{R}^{n\times q}$ is a boolean matrix that characterizes the placement of control devices and $\boldsymbol{u}\in\mathbb{R}^q$ denotes control input produced by the control devices; $t_0$ is an initial time and $\boldsymbol{x}_0$ and $\dot{\boldsymbol{x}}_0$ denote the initial displacement and velocity.
\end{sloppypar}

The above second-order system can be transformed into a first-order system as follows
\begin{align}
\label{eq:full-first}
&\boldsymbol{B}\dot{\boldsymbol{z}}	=\boldsymbol{A}\boldsymbol{z}+\boldsymbol{F}(\boldsymbol{z})+\epsilon\left(\boldsymbol{F}^\mathrm{ext}(t)+\boldsymbol{B}^{\mathrm{ext}}\boldsymbol{u}(t)\right),\nonumber \\
&\boldsymbol{z}(t_0)=\boldsymbol{z}_0,
\end{align}
where
\begin{align}
& \boldsymbol{z}=\begin{pmatrix}\boldsymbol{x}\\\dot{\boldsymbol{x}}\end{pmatrix},\,\,
\boldsymbol{A}=\begin{pmatrix}-\boldsymbol{K} 
& \boldsymbol{0}\\\boldsymbol{0} & \boldsymbol{M}\end{pmatrix}, \nonumber \\
& \boldsymbol{B}=\begin{pmatrix}\boldsymbol{C}_\mathrm{d} 
& \boldsymbol{M}\\\boldsymbol{M} & \boldsymbol{0}\end{pmatrix},\,\,
\boldsymbol{F}(\boldsymbol{z})=\begin{pmatrix}-\boldsymbol{f}(\boldsymbol{x},\dot{\boldsymbol{x}})\\\boldsymbol{0}\end{pmatrix},\nonumber \\
& \boldsymbol{F}^{\mathrm{ext}}(t)= \begin{pmatrix}\boldsymbol{E}(t)\\\boldsymbol{0}\end{pmatrix},\,\, \boldsymbol{B}^{\mathrm{ext}}= \begin{pmatrix}\boldsymbol{D}\\\boldsymbol{0}\end{pmatrix}.
\end{align}

Solving the linear part of~\eqref{eq:full-first} with $\epsilon=0$ leads to the generalized eigenvalue problem below
\begin{equation}
\boldsymbol{A}\boldsymbol{v}_j=\lambda_j\boldsymbol{B}\boldsymbol{v}_j,\quad \boldsymbol{u}_j^\ast \boldsymbol{A}=\lambda_j \boldsymbol{u}_j^\ast \boldsymbol{B},
\end{equation}
where $\lambda_j$ is a generalized eigenvalue and $\boldsymbol{v}_j$ and $\boldsymbol{u}_j$ are the corresponding \emph{right} and \emph{left} eigenvectors, respectively. Here we assume that the $\mathrm{Re}(\lambda_i)\neq0$ for $i=1,\cdots,2n$ such that the origin $\boldsymbol{z}=\boldsymbol{0}$ is a hyperbolic fixed point of the unforced system $\boldsymbol{B}\dot{\boldsymbol{z}}=\boldsymbol{A}\boldsymbol{z}+\boldsymbol{F}(\boldsymbol{z})$, namely, $\epsilon=0$.

\subsection{Vibration suppression as optimal control}
We aim to suppress the flexible structure's large amplitude nonlinear vibration via active control. In particular, we need to design the control input $\boldsymbol{u}(t)$ such that the response $\boldsymbol{z}(t)$ will be damped out quickly. Here, we use optimal control to achieve this goal. We define the objective functional below
\begin{equation}
\label{eq:obj-full}
J(\boldsymbol{u})=\int_{t_0}^{t_1}\left(\boldsymbol{z}^T\boldsymbol{Q}\boldsymbol{z}+\boldsymbol{u}^T\hat{\boldsymbol{R}}\boldsymbol{u}\right)dt +\boldsymbol{z}^T(t_1)\hat{\boldsymbol{M}}\boldsymbol{z}(t_1),
\end{equation}
where $t_0$ and $t_1$ are the initial and final times, $\boldsymbol{Q}$ and $\hat{\boldsymbol{M}}$ are symmetric positive semidefinite weight matrices, and $\hat{\boldsymbol{R}}$ is a symmetric positive definite weight matrix. Thus, the objective functional is the summation of the running cost $L(t,\boldsymbol{z},\boldsymbol{u})=\boldsymbol{z}^T\boldsymbol{Q}\boldsymbol{z}+\boldsymbol{u}^T\hat{\boldsymbol{R}}\boldsymbol{u}$ and the terminal cost $\boldsymbol{z}^T(t_1)\hat{\boldsymbol{M}}\boldsymbol{z}(t_1)$. The optimal control problem is formulated as below:
\begin{problem}\label{P1}
\textbf{High-dimensional nonlinear optimal control problem (HNOCP)}:
Find the control input $\boldsymbol{u}(t)$ and state $\boldsymbol{z}(t)$ that satisfy the state equations~\eqref{eq:full-first} to minimize the objective functional~\eqref{eq:obj-full}.
\end{problem}

\subsection{Challenges of optimal control design}
Although no path and control constraints have been considered, solving the above optimal control problem~\ref{P1} in real-time is challenging because of the high dimensionality and nonlinearity. Specifically, the nonlinearity makes it difficult to analyze the controllability of the system and solve the optimal control problem effectively. Moreover, the curse of high dimensionality ($n\gg1$) poses a significant challenge for real-time control. Indeed, solving an optimal control problem, even for a linear one, is computationally intensive when the state space is of high dimensions.

\begin{figure*}[!ht]
\centering
    \begin{tikzpicture}[node distance=3cm,
        rec/.style={rectangle, rounded corners, thick,minimum width=2.5cm,minimum height=1cm, text centered, draw=black,fill=#1!10, opacity=1},
        arrow/.style={ ->,>=stealth, line width=10pt, very thick, black!80!black}]
        \node(HNOCP) [rec=gray] {HNOCP};
        \node(HLOCP) [rec=gray, right of=HNOCP] {HLOCP};
        \node(LLOCP) [rec=gray, right of=HLOCP] {LLOCP};
        \node(ExLQR) [rec=gray, right of=LLOCP, align=center] {Extended\\ LQR};
        
        \draw[arrow]{} (HNOCP.east) --(HLOCP.west);
        \draw[arrow]{} (HLOCP.east) --(LLOCP.west);
        \draw[arrow]{} (LLOCP.east) --(ExLQR.west);
    \end{tikzpicture}
\caption{An SSM-based control design framework for the active control of nonlinear flexible structures. Here, HNOCP stands for high-dimensional \emph{nonlinear} optimal control problem~\ref{P1}, HLOCP denotes high-dimensional \emph{linear} optimal control problem~\ref{P2}, LLOCP represents \emph{low}-dimensional linear optimal control problem~\ref{P3}, and LQR stands for a linear quadratic regulator that solves problem~\ref{P4}.}
\label{ControlFramework}
\end{figure*}
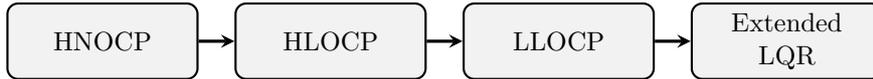

Next, we will show how to use spectral submanifolds (SSMs) to address the aforementioned two challenges. We propose an SSM-based control design framework that consists of three steps, as summarized in Fig.~\ref{ControlFramework}. In the first step, we decompose the state vector $\boldsymbol{z}$ into a nonlinear part and a linear part that corrects the nonlinear approximation. This decomposition enables us to transform the original \emph{nonlinear} optimal control problem~\ref{P1} into a \emph{linear} optimal control problem of the same dimensionality. In the second step, we perform a linear reduction to convert the \emph{high}-dimensional linear optimal control problem into a \emph{low}-dimensional linear optimal control problem. In the third step, we show that the reduced optimal control problem can be formulated as an extended linear quadratic problem, which can be solved by an extended linear quadratic regulator (LQR). Before we present details for these three steps, we introduce SSM-based model reduction in Sect.~\ref{sec:review-ssm}.

\section{SSM-based model reduction}
\label{sec:review-ssm}

In this section, we present essential ingredients for model reduction using spectral submanifolds (SSMs). For unforced systems, we have autonomous SSMs, which are unique, smoothest invariant manifolds tangent to some spectral subspaces under proper non-resonance conditions regarding the spectrum of the linear part of~\eqref{eq:full-first}, namely, $\boldsymbol{B}\dot{\boldsymbol{z}}=\boldsymbol{A}\boldsymbol{z}$~\cite{haller2016nonlinear}. In short, higher-dimensional SSMs are needed if internal resonances are presented within the spectrum of the linear system. In our study, only two or four-dimensional SSMs are involved. With the addition of the external forcing and control input, the autonomous SSMs are perturbed as temporally aperiodic SSMs~\cite{Haller2024NonlinearMR}, which provide a rigorous framework for constructing SSM-based reduced-order models (ROMs) with control.

\subsection{SSMs for unforced systems}
\label{sec:ssm-theory}
We consider the following $2m$-dimensional \emph{master} modal subspace 
\begin{equation}
\label{eq6}
\mathcal{E}=\Span\{\boldsymbol{v}^\mathcal{E}_1,\bar{\boldsymbol{v}}^\mathcal{E}_1,\cdots,\boldsymbol{v}^\mathcal{E}_m,\bar{\boldsymbol{v}}^\mathcal{E}_m\}.
\end{equation}
Here, we have assumed that the master subspace $\mathcal{E}$ is spanned by $m$-pairs of complex conjugate modes. Under the assumption of small damping, we have small real parts for the eigenvalues and near resonances between their imaginary parts. Specifically, we allow for the following type of near \emph{inner} or \emph{internal} resonances
\begin{equation}
\label{eq:res-inner}
\lambda_i^\mathcal{E}\approx\boldsymbol{l}\cdot\boldsymbol{\lambda}^\mathcal{E}+\boldsymbol{j}\cdot\bar{\boldsymbol{\lambda}}^\mathcal{E},\quad \bar{\lambda}_i^\mathcal{E}\approx\boldsymbol{j}\cdot\boldsymbol{\lambda}^\mathcal{E}+\boldsymbol{l}\cdot\bar{\boldsymbol{\lambda}}^\mathcal{E}
\end{equation}
for some $i\in\{1,\cdots,m\}$, where $\boldsymbol{l},\boldsymbol{j}\in\mathbb{N}_0^m,$ $\vert\boldsymbol{l}+\boldsymbol{j}\vert:=\sum_{k=1}^m (l_k+j_k)\geq2$, and
\begin{equation}
\boldsymbol{\lambda}^\mathcal{E}=(\lambda^\mathcal{E}_1,\cdots,\lambda^\mathcal{E}_m),\quad 
\bar{\boldsymbol{\lambda}}_\mathcal{E}=(\bar{\lambda}^\mathcal{E}_1,\cdots,\bar{\lambda}^\mathcal{E}_m).
\end{equation}
For example, for a mechanical system with 1:2 internal resonance between the first two pairs of modes, i.e., $\lambda_2^\mathcal{E}\approx2\lambda_1^\mathcal{E}$, we have $\lambda_2^\mathcal{E}\approx l_{21}\lambda_1^\mathcal{E}+l_{22}\lambda_2^\mathcal{E}+j_{21}\bar{\lambda}_1^\mathcal{E}+j_{22}\bar{\lambda}_2^\mathcal{E}$ for all $l_{ik},j_{ik}\in\mathbb{N}_0$ that satisfy $(l_{21}-j_{21})/2+l_{22}-j_{22}=1$.

We often select the master subspace $\mathcal{E}$ with eigenvalues of the largest real parts among the $2n$ eigenvalues. When all eigenvalues have negative real parts, we take a slow spectral subspace as the master subspace. In the case that some of the eigenvalues have positive real parts, we can take the unstable spectral subspace~\cite{cenedese2022data,liNonlinearModelReduction2023a} as the master subspace. One can also take a mix-mode spectral subspace whose spectrum has both positive and negative real parts to perform model reduction~\cite{haller2023nonlinear}.

The master subspace $\mathcal{E}$ is an invariant subspace to the linear system~$\boldsymbol{B}\dot{\boldsymbol{z}}=\boldsymbol{A}\boldsymbol{z}$. In other words, a trajectory initialized with a point in the subspace will remain in this subspace. Since $\mathcal{E}$ is attracting and invariant, a projection on $\mathcal{E}$ gives an exact model reduction for the linear system.

Under the addition of the nonlinear function $\boldsymbol{F}(\boldsymbol{z})$, the master subspace $\mathcal{E}$ is no longer invariant to the nonlinear system~$\boldsymbol{B}\dot{\boldsymbol{z}}=\boldsymbol{A}\boldsymbol{z}+\boldsymbol{F}(\boldsymbol{z})$. Instead, $\mathcal{E}$ is perturbed as some $2m$-dimensional invariant manifolds of the nonlinear system that are tangent to the master subspace. However, as illustrated in~\cite{haller2016nonlinear}, there are infinitely many such invariant manifolds for a given master subspace. Among them, there is a unique, smoothest invariant manifold, defined as the SSM associated with the master subspace and denoted by $\mathcal{W}(\mathcal{E})$.

The SSM can be viewed as an embedding of an open set in reduced coordinates $\boldsymbol{p}=(q_1,\bar{q}_1,\cdots,q_m,\bar{q}_m)$ via a map $\boldsymbol{z}=\boldsymbol{W}(\boldsymbol{p})$, where $q_i$ and $\bar{q}_i$ denote the parameterization coordinates corresponding to $\boldsymbol{v}_i^{\mathcal{E}}$ and $\bar{\boldsymbol{v}}_i^{\mathcal{E}}$, respectively. Moreover, the reduced dynamics on the SSM can be expressed as $\dot{\boldsymbol{p}}=\boldsymbol{R}(\boldsymbol{p})$. They satisfy the invariance equation below
\begin{equation}
\label{eq:invariance-auto}
 \boldsymbol{B}{D}_{\boldsymbol{p}}\boldsymbol{W}(\boldsymbol{p}) \boldsymbol{R}(\boldsymbol{p})
=\boldsymbol{A}\boldsymbol{W}(\boldsymbol{p})+\boldsymbol{F}(\boldsymbol{W}(\boldsymbol{p})).
\end{equation}
As an extension of a slow spectral master subspace, a slow SSM is an attracting invariant manifold and, hence, can be used to perform exact model reduction. Indeed, autonomous SSMs have enabled efficient and even analytic extraction of backbone curves and also prediction of free vibration of high-dimensional mechanical systems~\cite{jain2022compute,part-i,liModelReductionConstrained2023}.

\subsection{SSMs for generally non-autonomous systems}
\label{sec:time-dependent manifold}

Under further addition of the external forcing and control input $\epsilon(\boldsymbol{E}(t)+\boldsymbol{D}\boldsymbol{u}(t))$, the hyperbolic fixed point $\boldsymbol{z}=\boldsymbol{0}$ is perturbed into a unique nearby hyperbolic trajectory $\boldsymbol{z}^\ast(t)$ of the same stability type when the system is weakly forced, namely, $\epsilon\ll1$~\cite{Haller2024NonlinearMR}. Accordingly, the autonomous SSM $\mathcal{W}(\mathcal{E})$ is perturbed as an aperiodic SSM $\mathcal{W}(\mathcal{E}, t)$ that contains $\boldsymbol{z}(t)^\ast$ and acts as a locally invariant manifold of the non-autonomous system~\eqref{eq:full-first}~\cite{Haller2024NonlinearMR}.

Likewise, the invariant manifold $\mathcal{W}(\mathcal{E}, t)$ can be viewed as an embedding of an open set into the phase space of system~\eqref{eq:full-first} via a map $\boldsymbol{z}=\boldsymbol{W}_{\epsilon}(\boldsymbol{p},t)$, and the associated reduced dynamics on the aperiodic SSM is characterized by $\dot{\boldsymbol{p}} = \boldsymbol{R}_\epsilon(\boldsymbol{p},t)$. They satisfy the invariance equation below
\begin{align}
\boldsymbol{B}({D}_{\boldsymbol{p}} & \boldsymbol{W}_{\epsilon}(\boldsymbol{p},t) \boldsymbol{R}_\epsilon(\boldsymbol{p},t)+{D}_{t}\boldsymbol{W}_{\epsilon}(\boldsymbol{p},t)) \nonumber\\
& = \boldsymbol{A}\boldsymbol{W}_{\epsilon}(\boldsymbol{p},t)+\boldsymbol{F}(\boldsymbol{W}_{\epsilon}(\boldsymbol{p},t))\nonumber\\
& \quad + \epsilon\left(\boldsymbol{F}^\mathrm{ext}(t)+\boldsymbol{B}^{\mathrm{ext}}\boldsymbol{u}(t)\right).
\label{eq:invariance-tssm}
\end{align}
As we will see, the reduction on the aperiodic SSM enables us to construct lower-dimensional ROMs for control. Next, we present the computational methods of the SSMs, which lay the foundation for control design.

\subsection{Computation of SSMs}
\label{sec:comp-ssm}
Here, we do not follow the computational procedure proposed in~\cite{Haller2024NonlinearMR} but provide a simplified alternative to adapt our control design. We seek the unknown parameterization $\boldsymbol{W}_{\epsilon}(\boldsymbol{p},t)$ and vector field $\boldsymbol{R}_{\epsilon}(\boldsymbol{p},t)$ as an asymptotic series in $\epsilon$ given their smooth dependence on $\epsilon$:
\begin{gather}
    \boldsymbol{W}_{\epsilon}(\boldsymbol{p},t)=\boldsymbol{W}(\boldsymbol{p})+\epsilon \boldsymbol{X}(\boldsymbol{p},t)+\mathcal{O}(\epsilon^2),\label{eq:ssm-exp-eps}\\
    \boldsymbol{R}_{\epsilon}(\boldsymbol{p},t)=\boldsymbol{R}(\boldsymbol{p})+\epsilon \boldsymbol{S}(\boldsymbol{p},t)+\mathcal{O}(\epsilon^2)\label{eq:red-exp-eps}.
\end{gather}
Substituting the above expansions into the invariance equation~\eqref{eq:invariance-tssm} and collecting terms at $\mathcal{O}(\epsilon^0)$ yields
\begin{equation}
    \boldsymbol{B}{D}_{\boldsymbol{p}}\boldsymbol{W}(\boldsymbol{p})\boldsymbol{R}(\boldsymbol{p})=\boldsymbol{A}\boldsymbol{W}(\boldsymbol{p})+\boldsymbol{F}(\boldsymbol{W}(\boldsymbol{p})),
\label{eq:SSM-auto-eq}
\end{equation}
which is exactly the same as the invariance equation for autonomous SSM, as seen in~\eqref{eq:invariance-auto}. This invariance equation can be solved via a parameterization method implemented in SSMTool~\cite{ssmtool21}, an open-source package that supports the automated computation of the map $\boldsymbol{W}(\boldsymbol{p})$ and the vector field $\boldsymbol{R}(\boldsymbol{p)}$ up to any orders of expansion.

Furthermore, we obtain the invariance equation at $\mathcal{O}(\epsilon)$:
\begin{align}
\label{eq:SSM-nonauto-eq}
    & \boldsymbol{B}{D}_{\boldsymbol{p}}\boldsymbol{W}(\boldsymbol{p})\boldsymbol{S}(\boldsymbol{p},t)+\boldsymbol{B}D_{\boldsymbol{p}}\boldsymbol{X}(\boldsymbol{p},t)\boldsymbol{R}(\boldsymbol{p})\nonumber\\
    & + \boldsymbol{B}D_{t}\boldsymbol{X}(\boldsymbol{p},t) 
    =\boldsymbol{A}\boldsymbol{X}(\boldsymbol{p},t) +D\boldsymbol{F}(\boldsymbol{W}(\boldsymbol{p}))\boldsymbol{X}(\boldsymbol{p},t)\nonumber\\
    & +\boldsymbol{F}^\mathrm{ext}(t)+\boldsymbol{B}^{\mathrm{ext}}\boldsymbol{u}(t).
\end{align}
With $\boldsymbol{W}(\boldsymbol{p})$ and $\boldsymbol{R}(\boldsymbol{p})$ at hand, we solve~\eqref{eq:SSM-nonauto-eq} to obtain $\boldsymbol{X}(\boldsymbol{p},t)$ and $\boldsymbol{S}(\boldsymbol{p},t)$. Likewise, a Taylor expansion in $\boldsymbol{p}$ is used to approximate $\boldsymbol{X}$ and $\boldsymbol{S}$. The expansion coefficients are not constant but functions of $t$ and hence time-varying. As such, they can be expressed as
\begin{equation}
    \boldsymbol{X}(\boldsymbol{p},t)=\sum_{\boldsymbol{k}} \boldsymbol{X}_{\boldsymbol{k}}(t)\boldsymbol{p}^{\boldsymbol{k}},\quad 
    \boldsymbol{S}(\boldsymbol{p},t)=\sum_{\boldsymbol{k}} \boldsymbol{S}_{\boldsymbol{k}}(t)\boldsymbol{p}^{\boldsymbol{k}}.
\end{equation}

Here we restrict to a leading-order approximation in $\boldsymbol{p}$, as did in~\cite{jain2022compute,part-i}. Specifically, we let
\begin{equation}
\label{eq:leading-nonauto}
    \boldsymbol{X}(\boldsymbol{p},t)\approx \boldsymbol{X}_{\boldsymbol{0}}(t),\quad
    \boldsymbol{S}(\boldsymbol{p},t)\approx\boldsymbol{S}_{\boldsymbol{0}}(t)
\end{equation}
and then the reduced dynamics are of the form
\begin{equation}
\label{eq:red-nonauto-lead}
    \dot{\boldsymbol{p}}=\boldsymbol{R}(\boldsymbol{p})+\epsilon\boldsymbol{S}_{\boldsymbol{0}}(t)+\mathcal{O(\epsilon|\boldsymbol{p}|)}.
\end{equation}
Next, we derive the solution to $\boldsymbol{X}_0(t)$ and $\boldsymbol{S}_0(t)$. Substituting the leading-order approximation~\eqref{eq:leading-nonauto} into~\eqref{eq:SSM-nonauto-eq} yields
\begin{align}
\boldsymbol{B}{D}_{\boldsymbol{p}}&\boldsymbol{W}(\boldsymbol{p})\boldsymbol{S}_0(t)+\boldsymbol{B}\dot{\boldsymbol{X}}_{\boldsymbol{0}}(t)\nonumber\\
&=\boldsymbol{A}\boldsymbol{X}_{\boldsymbol{0}}(t) +D\boldsymbol{F}(\boldsymbol{W}(\boldsymbol{p}))\boldsymbol{X}_{\boldsymbol{0}}(t)\nonumber\\
&+\boldsymbol{F}^\mathrm{ext}(t)+\boldsymbol{B}^{\mathrm{ext}}\boldsymbol{u}(t).
\end{align}
Collecting the terms that are independent of $\boldsymbol{p}$ gives
\begin{align}
\label{eq:inv-x0-s0}
\boldsymbol{B}\dot{\boldsymbol{X}}_{\boldsymbol{0}}(t)&=\boldsymbol{A}\boldsymbol{X}_{\boldsymbol{0}}(t) +\boldsymbol{F}^\mathrm{ext}(t)+\boldsymbol{B}^{\mathrm{ext}}\boldsymbol{u}(t)\nonumber\\
&-\boldsymbol{B}\boldsymbol{W}_{\boldsymbol{I}}\boldsymbol{S}_0(t),
\end{align}
where $\boldsymbol{W}_{\boldsymbol{I}}$ stands for the linear expansion coefficient matrix of the map $\boldsymbol{z}=\boldsymbol{W}(\boldsymbol{p})$.

We note that the above equations are under-determined because both $\boldsymbol{X}_0$ and $\boldsymbol{S}_0$ are unknowns. Therefore, we have the freedom to choose the non-autonomous part $\boldsymbol{S}_{\boldsymbol{0}}(t)$. Next, we use this freedom to transform the original HNOCP~\ref{P1} into a linear optimal control problem of the same dimensionality.

\begin{remark}
\label{rk1}
    The expansion order for $\boldsymbol{W}(\boldsymbol{p})$ and $\boldsymbol{R}(\boldsymbol{p})$ can be determined via a convergence study of backbone curves for two-dimensional SSMs and an invariance error measure for higher-dimensional SSMs~\cite{liModelReductionConstrained2023}. For systems with a limit cycle in steady state, one can also check the convergence of the limit cycle to determine the expansion order~\cite{liNonlinearModelReduction2023a}.
\end{remark}

\begin{remark}
    The truncation at leading order of non-autonomous SSM shown in~\eqref{eq:leading-nonauto} assumes that the magnitude of $\epsilon\boldsymbol{X}_{\boldsymbol{0}}$ is much smaller than the autonomous part $\boldsymbol{W}(\boldsymbol{p})$, which is generally true if the external forcing and control input is small. However, the SSM persists for larger forcing~\cite{Haller2024NonlinearMR}. Numerical examples show that our control design framework works well even if $\epsilon\boldsymbol{X}_{\boldsymbol{0}}$ is of comparable magnitude as that of $\boldsymbol{W}(\boldsymbol{p})$.
\end{remark}

\section{Transforming nonlinear control as linear control via decomposition}
\label{sec:Controltransformation}

Since the control input $\boldsymbol{u}(t)$ is not in harmonic form in general and $[t_0,t_1]$ is a finite horizon, we can choose $\boldsymbol{S}_{\boldsymbol{0}}(t)=\boldsymbol{0}$ in~\eqref{eq:inv-x0-s0} such that the following holds
\begin{equation}
\label{eq:state-space-dynamics}
\boldsymbol{B}\dot{\boldsymbol{X}}_{\boldsymbol{0}}(t)=\boldsymbol{A}\boldsymbol{X}_{\boldsymbol{0}}(t) +\boldsymbol{F}^\mathrm{ext}(t)+\boldsymbol{B}^{\mathrm{ext}}\boldsymbol{u}(t).
\end{equation}
Then, the nonlinear state equations~\eqref{eq:full-first} with control is reduced to the linear state equations above with control. Importantly, this allows us to apply well-developed linear control theory to the nonlinear control problem. Specifically, we first perform an offline computation to simulate the reduced dynamics below without control
\begin{equation}
\label{eq:red-p-without-u}
    \dot{\boldsymbol{p}}=\boldsymbol{R}(\boldsymbol{p}),\quad \boldsymbol{p}(t_0)=\boldsymbol{p}_0.
\end{equation}
We can then control the response of $\boldsymbol{X}_{\boldsymbol{0}}(t)$ to achieve a desired full response $\boldsymbol{z}(t)$ via the decomposition below (cf.~\eqref{eq:ssm-exp-eps} and~\eqref{eq:leading-nonauto})
\begin{equation}
\label{eq:decomp-zt}
\boldsymbol{z}(t)=\boldsymbol{W}(\boldsymbol{p}(t))+\epsilon\boldsymbol{X}_{\boldsymbol{0}}(t).
\end{equation}

For a given initial state $\boldsymbol{z}_0$, we first determine the initial condition $\boldsymbol{p}_0$ for the reduced dynamics by a projection, and then the initial state $\boldsymbol{X}_0(t_0)$. In particular, we have
\begin{equation}
\label{eq:p0}
    \boldsymbol{p}_0 = 
      (\boldsymbol{U}^\mathcal{E})^\ast\boldsymbol{B}\boldsymbol{z}_0,
\end{equation}
where we have used the fact that $\boldsymbol{z}_0\approx\boldsymbol{V}^\mathcal{E}\boldsymbol{p}_0$ and $ (\boldsymbol{U}^\mathcal{E})^\ast\boldsymbol{B}\boldsymbol{V}=\boldsymbol{I}$ in the projection. Here,
\begin{align}
    \boldsymbol{V}^\mathcal{E} &= (\boldsymbol{v}^\mathcal{E}_1,\bar{\boldsymbol{v}}^\mathcal{E}_1,\cdots,\boldsymbol{v}^\mathcal{E}_m,\bar{\boldsymbol{v}}^\mathcal{E}_m),\nonumber \\
    \boldsymbol{U}^\mathcal{E} &= (\boldsymbol{u}^\mathcal{E}_1,\bar{\boldsymbol{u}}^\mathcal{E}_1,\cdots,\boldsymbol{u}^\mathcal{E}_m,\bar{\boldsymbol{u}}^\mathcal{E}_m).
\end{align}
We then obtain $\boldsymbol{X}_0(t_0)=\boldsymbol{z}_0-\boldsymbol{W}(\boldsymbol{p}_0)$. With the above, we can transform HNOCP (\ref{P1}) into a linear optimal control problem below:

\begin{problem}\label{P2}
\textbf{High-dimensional linear optimal control problem (HLOCP)}:
Find the control input $\boldsymbol{u}(t)$ and the state vector $\boldsymbol{z}(t)$ decomposed by~\eqref{eq:decomp-zt} to minimize the objective functional~\eqref{eq:obj-full} while subject the state equations~\eqref{eq:state-space-dynamics} with initial condition $\boldsymbol{X}_0(t_0)=\boldsymbol{z}_0-\boldsymbol{W}(\boldsymbol{p}_0)$. Here, the time history for the reduced coordinates $\boldsymbol{p}(t)$ can be computed prior in an offline stage using a forward simulation of the reduced dynamics~\eqref{eq:red-p-without-u} with an initial condition given by~\eqref{eq:p0}.
\end{problem}

\section{Reduction of high-dimensional linear control problem}
\label{sec:ControlL}

\subsection{Reducing the dimensionality of HLOCP}
The HLOCP (\ref{P2}) is still a high-dimensional optimal control problem when $n\gg1$. This high dimensionality makes real-time control infeasible. To resolve the curse of high-dimensionality, we further perform a linear reduction to the high-dimensional linear system~\eqref{eq:state-space-dynamics}. For convenience, we project the full state space of $\boldsymbol{X}_0$ to some $l$-dimensional spectral subspace $\mathcal{L}$ of the linear dynamics $\boldsymbol{B}\dot{\boldsymbol{X}}_0=\boldsymbol{A}\boldsymbol{X}_0$. In particular, we have
\begin{align}
\label{eq:exp-X0}
& \boldsymbol{X}_0(t) = \sum_{i=1}^l\hat{\boldsymbol{v}}_i q_i 
 = \hat{\boldsymbol{V}}\boldsymbol{q}(t),\nonumber \\
& \hat{\boldsymbol{V}} =(\hat{\boldsymbol{v}}_1,\cdots, \hat{\boldsymbol{v}}_l),\quad l\ll n,
\end{align}
where $\hat{\boldsymbol{v}}_i$ is the right eigenvector of the eigenvalue problem $\boldsymbol{A}\hat{\boldsymbol{v}}_i=\hat{\lambda}_i\boldsymbol{B}\hat{\boldsymbol{v}}_i$, and $\mathcal{L}$ is spanned by the column vectors of $\hat{\boldsymbol{V}}$. We note that $\hat{\boldsymbol{V}}$ is not necessarily equal to the master spectral subspace $\mathcal{E}$. Here, we select $\hat{\boldsymbol{V}}$ using linear reduction techniques developed in the control community. More details about the selection will be provided later. With the approximation~\eqref{eq:exp-X0}, the decomposition~\eqref{eq:decomp-zt} is updated as
\begin{equation}
\label{eq:decomp-zt-qt}
    \boldsymbol{z}(t)=\boldsymbol{W}(\boldsymbol{p}(t))+\epsilon\hat{\boldsymbol{V}}\boldsymbol{q}(t).
\end{equation}

We substitute~\eqref{eq:exp-X0} into~\eqref{eq:state-space-dynamics}, yielding
\begin{equation}
\label{eq:subs-x0-vq}
\boldsymbol{B}\hat{\boldsymbol{V}}\dot{\boldsymbol{q}}=\boldsymbol{A}\hat{\boldsymbol{V}}\boldsymbol{q}+\boldsymbol{F}^\mathrm{ext}(t)+\boldsymbol{B}^\mathrm{ext}\boldsymbol{u}(t).
\end{equation}
Let $\hat{\boldsymbol{u}}_i$ be the left eigenvector corresponding to $\hat{\lambda}_i$ and normalized according to $\hat{\boldsymbol{u}}_i^\ast\boldsymbol{B}\hat{\boldsymbol{v}}_i=1$, the orthonormalization $\hat{\boldsymbol{u}}_i^\ast\boldsymbol{B}\hat{\boldsymbol{v}}_j=\delta_{ij}$ holds. Let $\hat{\boldsymbol{U}}=(\hat{\boldsymbol{u}}_1,\cdots,\hat{\boldsymbol{u}}_l)$, we pre-multiply both sides of the equation~\eqref{eq:subs-x0-vq} with $\hat{\boldsymbol{U}}^\ast$, yielding
\begin{equation}   
\label{eq:state-qt-dynamics}
\dot{\boldsymbol{q}}=\hat{\boldsymbol{\Lambda}}\boldsymbol{q}+\hat{\boldsymbol{U}}^\ast\left(\boldsymbol{F}^\mathrm{ext}(t)+\boldsymbol{B}^\mathrm{ext}\boldsymbol{u}(t)\right),
\end{equation}
where $\hat{\Lambda}=\mathrm{diag}(\hat{\lambda}_1,\hat{\lambda}_2,\cdots,\hat{\lambda}_l)$.
The initial conditions for $\boldsymbol{q}$ are determined below
\begin{align} 
\label{eq:q0}
\boldsymbol{q}_0:&=\boldsymbol{q}(t_0) = \hat{\boldsymbol{U}}^\ast\boldsymbol{B}\boldsymbol{X}_{\boldsymbol{0}}(t_0) \nonumber\\
&=\hat{\boldsymbol{U}}^\ast\boldsymbol{B}\left(\boldsymbol{z}_0-\boldsymbol{W}(\boldsymbol{p}_0)\right)/\epsilon.
\end{align}
With the above, we now transform the HLOCP (\ref{P2}) into a low-dimensional optimal control problem below.

\begin{problem}\label{P3}
\textbf{Low-dimensional linear optimal control problem (LLOCP)}:
Find the control input $\boldsymbol{u}(t)$ and the state vector $\boldsymbol{z}(t)$ decomposed by~\eqref{eq:decomp-zt-qt} to minimize the objective functional~\eqref{eq:obj-full} while subject the state equations~\eqref{eq:state-qt-dynamics} with the initial condition~\eqref{eq:q0}. Here, the time history for the reduced coordinates $\boldsymbol{p}(t)$ can be computed prior in an offline stage using a forward simulation of the reduced dynamics~\eqref{eq:red-p-without-u} with an initial condition given by~\eqref{eq:p0}.
\end{problem}

\subsection{Selection of reduction basis}
\label{sec:selection-basis}
We are left with the selection of the projection basis $\hat{\boldsymbol{V}}$ shown in~\eqref{eq:exp-X0}. The accuracy of the linear reduction highly depends on the choice of $\hat{\boldsymbol{V}}$. In particular, we need to quantify the accuracy of the approximated linear system~\eqref{eq:state-qt-dynamics} and determine which eigenvector should be selected. Here, we present two selection schemes: one is based on direct current gains (DCgains), and the other is based on modal Hankel singular values (MHSVs).

\subsubsection{Selection based on DCgains}
We first present the transfer function between the control input and output of the linear system~\eqref{eq:state-space-dynamics}. We let $\boldsymbol{F}^\mathrm{ext}=\boldsymbol{0}$ and define an observation matrix $\boldsymbol{C}$ such that $\boldsymbol{y}=\boldsymbol{C}\boldsymbol{X}_0$, where $\boldsymbol{y}$ denotes the vector of observable of interest. Then the transfer function is obtained as below
\begin{equation}
    \boldsymbol{G}(s)=\boldsymbol{C}(s\boldsymbol{B}-\boldsymbol{A})^{-1}\boldsymbol{B}^\mathrm{ext}.
\end{equation}
Let $\boldsymbol{V}$ be the right eigenmatrix such that $\boldsymbol{A}\boldsymbol{V}=\boldsymbol{B}\boldsymbol{V}\boldsymbol{\Lambda}$, where $\boldsymbol{\Lambda}$ is a diagonal matrix composed of all eigenvalues. We denote the corresponding left eigenmatrix as $\boldsymbol{U}$ such that $\boldsymbol{U}^\ast\boldsymbol{B}\boldsymbol{V}=\boldsymbol{I}$ and $\boldsymbol{U}^\ast\boldsymbol{A}\boldsymbol{V}=\boldsymbol{\Lambda}$. Then the transfer function $\boldsymbol{G}(s)$ can be rewritten as below
\begin{align}
\label{eq:Gs}
     & \boldsymbol{G}(s)=\boldsymbol{C}\boldsymbol{V}(s\boldsymbol{I}-\boldsymbol{\Lambda})^{-1}\boldsymbol{U}^\ast\boldsymbol{B}^\mathrm{ext} = \sum_{i=1}^{2n}\boldsymbol{G}_i(s),\nonumber \\ 
     & \boldsymbol{G}_i(s)=\frac{(\boldsymbol{C}\boldsymbol{v}_i) (\boldsymbol{u}_i^\ast\boldsymbol{B}^\mathrm{ext})}{s-\lambda_i},
\end{align}
where $\boldsymbol{G}_i(s)$ represents the contribution of the $i$-th eigenmode, $\boldsymbol{v}_i$ and $\boldsymbol{u}_i$ are the $i$-th column of the matrices of $\boldsymbol{V}$ and $\boldsymbol{U}$, and $\lambda_i$ is the associated eigenvalue. 

We observe from the decoupling in~\eqref{eq:Gs} that we can use the transfer function to characterize the contribution of each eigenmode. Following this idea, DCgain that depends on $\boldsymbol{G}(0)$ characterizes the steady-state ratio between input and output and can be used to measure the importance of each eigenmode. DCgain is also directly related to the infinity-induced norm of the system as a notion of output reachability~\cite{rantzerDistributedControlPositive2011, kawanoDataDrivenModelReduction2020}. Thus, it has been used to truncate not-so-important state variables and perform model reduction~\cite{kawanoDataDrivenModelReduction2020}.

In our setting, the eigenvalues appear in complex conjugate pairs. It follows that~\eqref{eq:Gs} can be further rewritten as
\begin{align}
\label{eq:Gs-n}
     & \boldsymbol{G}(s)= \sum_{i=1}^{n}\hat{\boldsymbol{G}}_i(s),\nonumber \\
     & \hat{\boldsymbol{G}}_i(s) =\frac{(\boldsymbol{C}\boldsymbol{v}_i) (\boldsymbol{u}_i^\ast\boldsymbol{B}^\mathrm{ext})}{s-\lambda_i}+\frac{(\boldsymbol{C}\bar{\boldsymbol{v}}_i) (\bar{\boldsymbol{u}}_i^\ast\boldsymbol{B}^\mathrm{ext})}{s-\bar{\lambda}_i}.
\end{align}
Then, the contribution of the $i$-th pair of complex conjugate eigenmodes to the DCgain is computed below
\begin{equation}
    \text{DCgain}_i=\lim_{s\to 0}||\hat{\boldsymbol{G}}_i(s)||_2.
\end{equation}
and its associated significance can be characterized by
\begin{equation}
{p}_i=\frac{\mathrm{DCgain}_i}{\sum_{j=1}^n\mathrm{DCgain}_j}.
\end{equation}
For high-dimensional systems, it is infeasible to compute the DCgains for all eigenmodes. We often use a truncation up to the first $\hat{m}$ pairs of modes with the lowest frequencies when the sequence of DCgains $\{p_j\}$ is nearly converged; namely,
\begin{equation}
{p}_i=\frac{\mathrm{DCgain}_i}{\sum_{j=1}^{\hat{m}}\mathrm{DCgain}_j}.
\end{equation}
We refer $p_i$ above as the normalized DCgain of the $i$-th pair of modes. Indeed, high-frequency components of linear systems are difficult to excite in most operating conditions, which suggests the contribution of high-frequency eigenvectors can often be ignored. Then, we can select the $l$ pairs of modes of largest contributions such that the summation of their significance is larger than a threshold that quantifies the accuracy of the reduction using  the $l$ selected pairs of modes.

\subsubsection{Selection based on MHSVs}

The selection based on DCgains uses $\hat{\boldsymbol{G}}_i(0)$, namely, the transfer function evaluated at zero. An alternative approach is to use the transfer function $\hat{\boldsymbol{G}}_i(s)$ for whole $s\in\mathbb{C}$. Following this idea, modal Hankel singular values (MHSVs) have been used to quantify the significance of eigenmodes~\cite{Chang2002DesignOR,changModelReductionBased2004}. 

Let $\lambda_i$ and $\bar{\lambda}_i$ be the eigenvalues of a pair of complex conjugate eigenmodes, the associated modal coordinates, right and left eigenvectors are denoted as $\tilde{\boldsymbol{q}}_i=(q_i,\bar{q}_i)$, $\tilde{\boldsymbol{V}}_i=(\boldsymbol{v}_i,\bar{\boldsymbol{v}}_i)$, and $\tilde{\boldsymbol{U}}_i=(\boldsymbol{u}_i,\bar{\boldsymbol{u}}_i)$. The equations of motion for $\tilde{\boldsymbol{q}}_i$ are listed below
\begin{equation}   
\label{eq:state-qt-dynamics-individual}
\dot{\tilde{\boldsymbol{q}}}_i=\tilde{\boldsymbol{\Lambda}}_i\tilde{\boldsymbol{q}}_i+\tilde{\mathbf{B}}_i\boldsymbol{u}(t),
\end{equation}
where $\tilde{\boldsymbol{\Lambda}}_i=\mathrm{diag}\{{\lambda}_i,\bar{{\lambda}}_i \}$ and $\tilde{\mathbf{B}}_i=\tilde{\boldsymbol{U}}_i^\ast\boldsymbol{B}^\mathrm{ext}$. The corresponding observation matrix is obtained as $\tilde{\boldsymbol{C}}_i=\boldsymbol{C}\tilde{\boldsymbol{V}}_i$ because $\boldsymbol{y}=\boldsymbol{C}\boldsymbol{X}_0\approx\boldsymbol{C}\tilde{\boldsymbol{V}}_i\tilde{\boldsymbol{q}}_i=\tilde{\boldsymbol{C}}_i\tilde{\boldsymbol{q}}_i$. Provided that $\mathrm{Re}(\lambda_i)<0$, the controllability Gramian and observability Gramian of the system $(\tilde{\boldsymbol{\Lambda}}_i,\tilde{\boldsymbol{B}}_i,\tilde{\boldsymbol{C}}_i)$ are as below:
\begin{align}
\label{eq:WcandWo}    \boldsymbol{W}_i^\mathrm{C}&=\int_{t_0}^{\infty}\boldsymbol{e}^{\tilde{\boldsymbol{\Lambda}}_it}\tilde{\boldsymbol{B}}_i\tilde{\boldsymbol{B}}_i^T\boldsymbol{e}^{\tilde{\boldsymbol{\Lambda}}_i^Tt}dt, \nonumber \\
\boldsymbol{W}_i^\mathrm{O}&=\int_{t_0}^{\infty}\boldsymbol{e}^{\tilde{\boldsymbol{\Lambda}}_i^Tt}\tilde{\boldsymbol{C}}_i^T\tilde{\boldsymbol{C}}_i\boldsymbol{e}^{\tilde{\boldsymbol{\Lambda}}_it}dt.
\end{align}
They are given as the solutions of the following Lyapunov equations
\begin{align}
\tilde{\boldsymbol{\Lambda}}_i\boldsymbol{W}_i^\mathrm{C}+\boldsymbol{W}_i^\mathrm{C}\tilde{\boldsymbol{\Lambda}}_i^T+\tilde{\boldsymbol{B}}_i\tilde{\boldsymbol{B}}_i^T &=\boldsymbol{0},\nonumber \\
\tilde{\boldsymbol{\Lambda}}_i^T\boldsymbol{W}_i^\mathrm{O}+\boldsymbol{W}_i^\mathrm{O}\tilde{\boldsymbol{\Lambda}}_i+\tilde{\boldsymbol{C}}_i\tilde{\boldsymbol{C}}_i^T &=\boldsymbol{0}.
\end{align}
Then the balanced Gramian $\mathrm{diag}\{\sigma_{i}^2,(\sigma_{i}')^2\}$ of the $i$-th pair of modes is defined as the Hankel singular values of the product of both Gramians~\cite{Kung1980OptimalHM,Moore1981PrincipalCA,changModelReductionBased2004}
\begin{equation}   \boldsymbol{W}_i^\mathrm{C}\cdot\boldsymbol{W}_i^\mathrm{O}=\check{\boldsymbol{U}}_i\begin{pmatrix}\sigma_i^2 & 0\\0 & (\sigma_i')^2\end{pmatrix}\check{\boldsymbol{V}}_i,\quad 
\end{equation}
and the MHSV of the $i$-th pair of modes is given by $\sigma_i^M=\max(\sigma_i,\sigma_i')$.

If we keep the first $\hat{m}$ pairs of eigenmodes, the model approximation error bound is given by~\cite{Kung1980OptimalHM,Chang2002DesignOR}
\begin{equation}
    ||\boldsymbol{G}(s)-\hat{\boldsymbol{G}}(s)||_{\infty}\leq 4\sum_{i=\hat{m}+1}^{n}\sigma_i^M, 
\end{equation}
where $\hat{\boldsymbol{G}}(s)$ is the truncated approximation of the original transfer function ${\boldsymbol{G}}(s)$ via the first $\hat{m}$ pairs of modes and $||\cdot||_{\infty}$ is the $\mathcal{H}_{\infty}$ norm. Therefore, the MHSV $\sigma_i^M$ characterizes the significance of the $i$-th pair of modes. Similar to the previous case, we use a truncation up to the first $\hat{m}$ ($\hat{m}\ll n$) pairs of modes with the lowest frequencies for approximating $\boldsymbol{G}(s)$ of high-dimensional systems. We calculate the normalized MHSV of each of these modes and then select $l$ ($l\leq\hat{m}$) pairs of modes of largest contributions such that the summation of their significance is larger than a threshold.

\begin{remark}
\label{rk:mhsv}
    DCgain and MSHV are well defined only for stable modes, as indicated by~\eqref{eq:WcandWo}, where the integral is undefined when the real part of $\lambda_i$ is positive. Here, we include all unstable eigenmodes in the reduction basis and compute DCgain and MHSV for stable eigenmodes. We further select a subset of stable eigenmodes into the reduction basis using the computed DCgain and MHSV. This is reasonable because only stable modes decay quickly.
\end{remark}

\section{Solutions to reduced optimal control problem}
\label{sec6:LQ}

\subsection{Linear quadratic (LQ) formulation}
We first show that the LLOCP (\ref{P3}) is actually an extended LQ problem. We recall that the running cost in the objective functional~\eqref{eq:obj-full} is given as $L(t,\boldsymbol{z},\boldsymbol{u})=\boldsymbol{z}^T\boldsymbol{Q}\boldsymbol{z}+\boldsymbol{u}^T\hat{\boldsymbol{R}}\boldsymbol{u}$. With the decomposition~\eqref{eq:decomp-zt-qt}, we have
\begin{align}
\boldsymbol{z}^T\boldsymbol{Q}\boldsymbol{z} & =\left(\boldsymbol{W}(\boldsymbol{p})+\epsilon\hat{\boldsymbol{V}}\boldsymbol{q}\right)^T\boldsymbol{Q}\left(\boldsymbol{W}(\boldsymbol{p})+\epsilon\hat{\boldsymbol{V}}\boldsymbol{q}\right)\nonumber\\
& = (\boldsymbol{W}(\boldsymbol{p}))^T\boldsymbol{Q}\boldsymbol{W}(\boldsymbol{p})+2\epsilon (\boldsymbol{W}(\boldsymbol{p}))^T\boldsymbol{Q}\hat{\boldsymbol{V}}\boldsymbol{q}\nonumber\\
& +\epsilon^2\boldsymbol{q}^T\hat{\boldsymbol{V}}^T\boldsymbol{Q}\hat{\boldsymbol{V}}\boldsymbol{q}\nonumber\\
& = a(t)+\boldsymbol{b}_{\boldsymbol{Q}}^T(t)\boldsymbol{q}+\boldsymbol{q}^T\boldsymbol{Q}_2\boldsymbol{q},
\end{align}
where
\begin{align}
\label{eq:abQ}
&a(t)=(\boldsymbol{W}(\boldsymbol{p}))^T\boldsymbol{Q}\boldsymbol{W}(\boldsymbol{p}),\nonumber \\
&\boldsymbol{b}_{\boldsymbol{Q}}(t)=2\epsilon (\boldsymbol{Q}\hat{\boldsymbol{V}})^T\boldsymbol{W}(\boldsymbol{p}),\nonumber \\
&\boldsymbol{Q}_2=\epsilon^2\hat{\boldsymbol{V}}^T\boldsymbol{Q}\hat{\boldsymbol{V}}.
\end{align}
Likewise, the terminal cost can be rewritten as
\begin{align}
\boldsymbol{z}^T(t_1)\hat{\boldsymbol{M}}\boldsymbol{z}(t_1)&= a(t_1) +\boldsymbol{b}_{\boldsymbol{M}}^T(t_1)\boldsymbol{q}(t_1)\nonumber\\
& +\boldsymbol{q}^T(t_1)\boldsymbol{M}_2\boldsymbol{q}(t_1)
\end{align}
where
\begin{equation}
\label{eq:bM}
\boldsymbol{b}_{\boldsymbol{M}}(t)=2\epsilon (\hat{\boldsymbol{M}}\hat{\boldsymbol{V}})^T\boldsymbol{W}(\boldsymbol{p}),\quad \boldsymbol{M}_2=\epsilon^2\hat{\boldsymbol{V}}^T\hat{\boldsymbol{M}}\hat{\boldsymbol{V}}.
\end{equation}

We recall that the solution $\boldsymbol{p}(t)$ does not depend on the control input $\boldsymbol{u}$ and can be computed in priori. Therefore, we can think of ${a}(t)$, $\boldsymbol{b}_{\boldsymbol{Q}}(t)$, and $\boldsymbol{b}_{\boldsymbol{M}}(t)$ as known signals. Consequently, the objective functional $J$ in~\eqref{eq:obj-full} is equivalent to $\tilde{J}$ below
\begin{align}
\label{eq:Jtilde}
\tilde{J}(\boldsymbol{u}) = & \int_{t_0}^{t_1} \left(\boldsymbol{b}_{\boldsymbol{Q}}^T(t)\boldsymbol{q}+\boldsymbol{q}^T\boldsymbol{Q}_2\boldsymbol{q}+\boldsymbol{u}^T\hat{\boldsymbol{R}}\boldsymbol{u}\right)dt \nonumber \\
& +\boldsymbol{b}_{\boldsymbol{M}}^T(t_1)\boldsymbol{q}(t_1)+\boldsymbol{q}^T(t_1)\boldsymbol{M}_2\boldsymbol{q}(t_1).
\end{align}
Here, we have dropped out the constant terms related to $a(t)$ because it has no effect on the optimal control law. Now, we are ready to present the extended LQ problem formulation below
\begin{problem}\label{P4}
\textbf{Extended linear quadratic problem (ELQP)}:
Find the control input $\boldsymbol{u}(t)$ and the state vector ${\boldsymbol{q}}(t)$ to minimize the quadratic objective functional~\eqref{eq:Jtilde} while subject the linear state equations~\eqref{eq:state-qt-dynamics} with the initial condition~\eqref{eq:q0}.
\end{problem}

\subsection{Solution to ELQP as an extended linear–quadratic regulator (LQR)}

We use calculus of variation to solve the optimal control problem~\ref{P4}. Let
\begin{gather}
\hat{\boldsymbol{B}} = \hat{\boldsymbol{U}}^\ast\boldsymbol{B}^\mathrm{ext},\quad \boldsymbol{b}(t)=\hat{\boldsymbol{U}}^\ast\boldsymbol{F}^\mathrm{ext}(t),\nonumber\\
\hat{L}(t,\boldsymbol{q},\boldsymbol{u})=\boldsymbol{b}_{\boldsymbol{Q}}^T(t)\boldsymbol{q}+\boldsymbol{q}^T\boldsymbol{Q}_2\boldsymbol{q}+\boldsymbol{u}^T\hat{\boldsymbol{R}}\boldsymbol{u},\nonumber\\
\Phi(\boldsymbol{q}(t_1))=\boldsymbol{b}_{\boldsymbol{M}}^T(t_1)\boldsymbol{q}(t_1)+\boldsymbol{q}^T(t_1)\boldsymbol{M}_2\boldsymbol{q}(t_1),
\end{gather}
we define an augmented Lagrangian below
\begin{align}
\hat{J} & = \int_{t_0}^{t_1} \left(\hat{L}(t,\boldsymbol{q},\boldsymbol{u})\right.\nonumber \\ & 
\left.\quad+ \boldsymbol{\mu}^\ast\left(\dot{\boldsymbol{q}}-\hat{\boldsymbol{\Lambda}}\boldsymbol{q}-\hat{\boldsymbol{B}}\boldsymbol{u}-\boldsymbol{b}(t)\right)\right) dt +\Phi(\boldsymbol{q}(t_1))\nonumber\\ 
& = \int_{t_0}^{t_1} \left( \boldsymbol{\mu}^\ast\dot{\boldsymbol{q}}-H(t,\boldsymbol{q},\boldsymbol{u},\bar{\boldsymbol{\mu}})\right)dt +\Phi(\boldsymbol{q}(t_1)),
\end{align}
where the Hamiltonian is given by
\begin{equation}
H(t,\boldsymbol{q},\boldsymbol{u},\bar{\boldsymbol{\mu}}) =\bar{\boldsymbol{\mu}}^T\left(\hat{\boldsymbol{\Lambda}}\boldsymbol{q}+\hat{\boldsymbol{B}}\boldsymbol{u}+\boldsymbol{b}(t)\right)-\hat{L}(t,\boldsymbol{q},\boldsymbol{u}).
\end{equation}

Calculus of variation yields
\begin{align}
\delta\hat{J} = & \int_{t_0}^{t_1} \left(\delta\boldsymbol{\mu}^\ast\dot{\boldsymbol{q}}+\boldsymbol{\mu}^\ast \delta\dot{\boldsymbol{q}}- H_{\boldsymbol{q}}^T\delta\boldsymbol{q}
\right.\nonumber\\
& \left.-H_{\bar{\boldsymbol{\mu}}}^T\delta\bar{\boldsymbol{\mu}}- H_{\boldsymbol{u}}^T\delta\boldsymbol{u}\right)dt+\Phi_{\boldsymbol{q}(t_1)}^T\delta\boldsymbol{q}(t_1)\nonumber\\
= & \int_{t_0}^{t_1} \left(\delta\boldsymbol{\mu}^\ast\dot{\boldsymbol{q}}-\dot{\boldsymbol{\mu}}^\ast \delta{\boldsymbol{q}}-\delta\boldsymbol{q}^TH_{\boldsymbol{q}} -\delta\bar{\boldsymbol{\mu}}^TH_{\bar{\boldsymbol{\mu}}}\right.\nonumber\\
& \left.-\delta\boldsymbol{u}^TH_{\boldsymbol{u}}\right)dt+(\bar{\boldsymbol{\mu}}^T(t_1)+\Phi_{\boldsymbol{q}(t_1)}^T)\delta\boldsymbol{q}(t_1),
\end{align}
where $H_{\boldsymbol{a}}$ gives the gradient of $H$ with respect to $\boldsymbol{a}$ for $\boldsymbol{a}=\{\boldsymbol{q},\bar{\boldsymbol{\mu}},\boldsymbol{u}\}$. Vanishing $\delta\hat{J}$ yields a set of first-order necessary conditions:
\begin{align}
& \delta\bar{\boldsymbol{\mu}}: \dot{\boldsymbol{q}}=H_{\boldsymbol{\mu}}=\hat{\boldsymbol{\Lambda}}\boldsymbol{q}+\hat{\boldsymbol{B}}\boldsymbol{u}+\boldsymbol{b}(t),\label{eq:Hmu}\\
& \delta{\boldsymbol{q}}: \dot{\bar{\boldsymbol{\mu}}}=-H_{\boldsymbol{q}}=\boldsymbol{b}_{\boldsymbol{Q}}(t)+2\boldsymbol{Q}_2\boldsymbol{q}-\hat{\boldsymbol{\Lambda}}\bar{\boldsymbol{\mu}},\nonumber \\ & \delta{\boldsymbol{q}}(t_1):  \bar{\boldsymbol{\mu}}(t_1)=-\Phi_{\boldsymbol{q}(t_1)},\label{eq:Hq}\\
&\delta\boldsymbol{u}: H_{\boldsymbol{u}}=-2\hat{\boldsymbol{R}}\boldsymbol{u}+\hat{\boldsymbol{B}}^T\bar{\boldsymbol{\mu}}=\boldsymbol{0}.\label{eq:Hu}
\end{align}
We solve the control input from~\eqref{eq:Hu} that
\begin{equation}
\label{eq:u-mu}
\boldsymbol{u}=\frac{1}{2}\hat{\boldsymbol{R}}^{-1}\hat{\boldsymbol{B}}^T\bar{\boldsymbol{\mu}}.
\end{equation}
We substitute the equation above into~\eqref{eq:Hmu} and~\eqref{eq:Hq}, yielding
\begin{gather}
\dot{\boldsymbol{q}} = \hat{\boldsymbol{\Lambda}}\boldsymbol{q}+\frac{1}{2}\hat{\boldsymbol{B}}\hat{\boldsymbol{R}}^{-1}\hat{\boldsymbol{B}}^T\bar{\boldsymbol{\mu}}+\boldsymbol{b}(t),\quad \boldsymbol{q}(t_0)=\boldsymbol{q}_0,\label{eq:dotq}\\
\dot{\bar{\boldsymbol{\mu}}}=2\boldsymbol{Q}_2\boldsymbol{q}-\hat{\boldsymbol{\Lambda}}\bar{\boldsymbol{\mu}}+\boldsymbol{b}_{\boldsymbol{Q}}(t),\nonumber\\
\bar{\boldsymbol{\mu}}(t_1)=-\boldsymbol{b}_{\boldsymbol{M}}(t_1)-2\boldsymbol{M}_2\boldsymbol{q}(t_1).\label{eq:dotmu}
\end{gather}

To obtain a feedback control law, we set
\begin{equation}
\label{eq:mubar}
  \bar{\boldsymbol{\mu}}=-2\mathbf{P}(t)\boldsymbol{q}+\boldsymbol{s}(t)  
\end{equation}
where $\mathbf{P}\in\mathbb{C}^{l\times l}$ and $\boldsymbol{s}\in\mathbb{C}^l$ are unknowns. Compared to classic LQR, we add the augmentation vector $\boldsymbol{s}$ to extend the classic LQR to our feedback control design. Next, we show how to solve for the matrix $\mathbf{P}$ and the vector $\boldsymbol{s}$.

Differentiation of~\eqref{eq:mubar} gives
\begin{equation}
\dot{\bar{\boldsymbol{\mu}}} = -2\dot{\mathbf{P}}\boldsymbol{q}-2\mathbf{P}\dot{\boldsymbol{q}}+\dot{\boldsymbol{s}}.
\end{equation}
Substituting the canonical equations~\eqref{eq:dotq} and~\eqref{eq:dotmu} into the above equation yields
\begin{align}
& 2\boldsymbol{Q}_2\boldsymbol{q}-\hat{\boldsymbol{\Lambda}}(-2\mathbf{P}\boldsymbol{q}+\boldsymbol{s})+\boldsymbol{b}_{\boldsymbol{Q}}=-2\dot{\mathbf{P}}\boldsymbol{q} \nonumber \\
& -2\mathbf{P}\left(\hat{\boldsymbol{\Lambda}}\boldsymbol{q}+\frac{1}{2}\hat{\boldsymbol{B}}\hat{\boldsymbol{R}}^{-1}\hat{\boldsymbol{B}}^T(-2\mathbf{P}\boldsymbol{q}+\boldsymbol{s})+\boldsymbol{b}\right)+\dot{\boldsymbol{s}}.
\end{align}
We collect the terms at $\mathcal{O}(1)$ and $\mathcal{O}(\boldsymbol{q})$ in the above equation, yielding
\begin{align}
\mathcal{O}(1): &-\hat{\boldsymbol{\Lambda}}\boldsymbol{s}+\boldsymbol{b}_{\boldsymbol{Q}}=-\mathbf{P}\hat{\boldsymbol{B}}\hat{\boldsymbol{R}}^{-1}\hat{\boldsymbol{B}}^T\boldsymbol{s}-2\mathbf{P}\boldsymbol{b}+\dot{\boldsymbol{s}},\label{eq:o1}\\
\mathcal{O}(\boldsymbol{q}): & 2\boldsymbol{Q}_2+2\hat{\boldsymbol{\Lambda}}\mathbf{P}=-2\dot{\mathbf{P}}-2\mathbf{P}\hat{\boldsymbol{\Lambda}}\nonumber\\
&\quad +2\mathbf{P}\hat{\boldsymbol{B}}\hat{\boldsymbol{R}}^{-1}\hat{\boldsymbol{B}}^T\mathbf{P}\label{eq:oq}.
\end{align}
Here,~\eqref{eq:oq} is the Riccati equation:
\begin{equation}
\label{eq:riccati}
\dot{\mathbf{P}} = -\mathbf{P}\hat{\boldsymbol{\Lambda}}-\hat{\boldsymbol{\Lambda}}\mathbf{P}-\boldsymbol{Q}_2+\mathbf{P}\hat{\boldsymbol{B}}\hat{\boldsymbol{R}}^{-1}\hat{\boldsymbol{B}}^T\mathbf{P},
\end{equation}
while~\eqref{eq:o1} gives a time-varying ODE for the augmentation vector~$\boldsymbol{s}$:
\begin{equation}
\label{eq:dots}
\dot{\boldsymbol{s}} = \left(\mathbf{P}\hat{\boldsymbol{B}}\hat{\boldsymbol{R}}^{-1}\hat{\boldsymbol{B}}^T-\hat{\boldsymbol{\Lambda}}\right)\boldsymbol{s}+\boldsymbol{b}_{\boldsymbol{Q}}+2\mathbf{P}\boldsymbol{b}.
\end{equation}
To solve for~\eqref{eq:riccati} and~\eqref{eq:dots}, we need initial or terminal conditions for $\mathbf{P}$ and $\boldsymbol{s}$. We deduce from
\begin{align}
 \bar{\boldsymbol{\mu}}(t_1)&=-\boldsymbol{b}_{\boldsymbol{M}}(t_1)-2\boldsymbol{M}_2\boldsymbol{q}(t_1)\nonumber\\
 &=-2\mathbf{P}(t_1)\boldsymbol{q}(t_1)+\boldsymbol{s}(t_1)
\end{align}
and $\boldsymbol{q}(t_1)$ is free to change that 
\begin{equation}
\label{eq:final-Pands}
\mathbf{P}(t_1) = \boldsymbol{M}_2,\quad \boldsymbol{s}(t_1)=-\boldsymbol{b}_{\boldsymbol{M}}(t_1).
\end{equation}
Then, we can use backward simulations to obtain $\mathbf{P}(t)$ and $\boldsymbol{s}(t)$.

Consequently, the control input is given as
\begin{equation}
\label{eq:u-closed-loop-q}
\boldsymbol{u}=-\hat{\boldsymbol{R}}^{-1}\hat{\boldsymbol{B}}^T\mathbf{P}(t)\boldsymbol{q}+0.5\hat{\boldsymbol{R}}^{-1}\hat{\boldsymbol{B}}^T\boldsymbol{s}(t).
\end{equation}
Moving back to physical coordinates, we obtain
\begin{align}
\boldsymbol{u}=&-\hat{\boldsymbol{R}}^{-1}\hat{\boldsymbol{B}}^T\mathbf{P}(t)\hat{\boldsymbol{U}}^\ast\boldsymbol{B}\left[(\boldsymbol{z}-\boldsymbol{W}(\boldsymbol{p}))\right]/\epsilon \nonumber\\
&+0.5\hat{\boldsymbol{R}}^{-1}\hat{\boldsymbol{B}}^T\boldsymbol{s}(t).
\end{align}

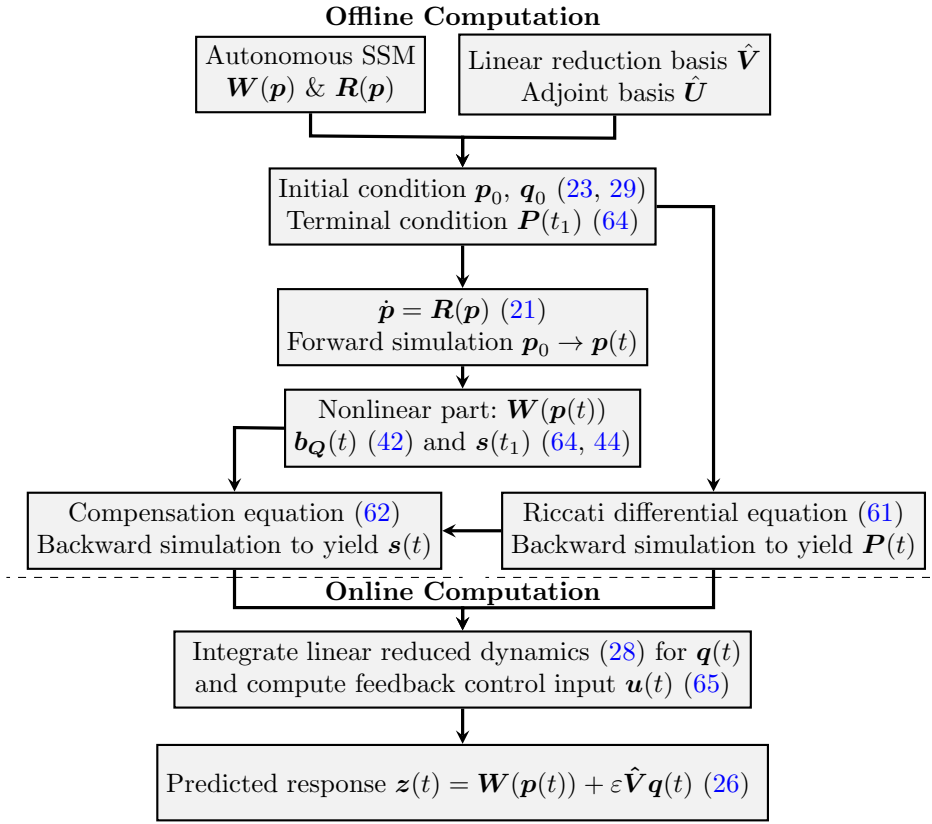
\begin{figure*}[!ht]
\centering
    \begin{tikzpicture}[node distance=3cm,
    rec/.style={rectangle, thick,minimum width=2.5cm,minimum height=1cm, text centered, draw=black,fill=#1!10, opacity=1},
    arrow/.style={ ->,>=stealth, very thick, black!80!black}]
        \node(note1) [yshift=0.75cm, xshift=0cm]{\textbf{Offline Computation}};
        \node(note2) [yshift=-6.875cm, xshift=0cm]{\textbf{Online Computation}};
        \node(SSM) [rec=gray, align=center, xshift=-2cm] {Autonomous SSM \\
        $\boldsymbol{W}(\boldsymbol{p})$ \& $\boldsymbol{R}(\boldsymbol{p})$ };
        \node(uv)  [rec=gray, right of=SSM, align=center, xshift=1cm]  {Linear reduction basis $\hat{\boldsymbol{V}}$ \\
        Adjoint basis $\hat{\boldsymbol{U}}$
        };
        \node(IC) [rec=gray, below of=SSM, xshift=2cm, yshift=1.25cm, align=center]  {
        Initial condition 
        $\boldsymbol{p}_0$, $\boldsymbol{q}_0$ (\ref{eq:p0}, \ref{eq:q0}) \\
        Terminal condition
        $\boldsymbol{P}(t_1)$ (\ref{eq:final-Pands})    
        };
        \node(calP) [rec=gray, below of=IC, yshift=1.4cm, align=center]  {
        $\dot{\boldsymbol{p}}=\boldsymbol{R}(\boldsymbol{p})$ (\ref{eq:red-p-without-u}) \\
        Forward simulation $\boldsymbol{p}_0\to \boldsymbol{p}(t)$
        };
        \node(WbQbM) [rec=gray, below of=calP, align=center, yshift=1.66cm] {
        Nonlinear part:~$ \boldsymbol{W}(\boldsymbol{p}(t))$ \\
        $\boldsymbol{b}_{\boldsymbol{Q}}(t)$ (\ref{eq:abQ}) and $\boldsymbol{s}(t_1)$ (\ref{eq:final-Pands}, \ref{eq:bM})
        };
        \node(Ricti) [rec=gray, below of=IC, xshift=3.3cm, yshift=-1.3cm, align=center] {
        Riccati differential equation (\ref{eq:riccati}) \\ 
        Backward simulation to yield $\boldsymbol{P}(t)$
        };
        \node(st) [rec=gray, below of=IC, yshift=-1.3cm, xshift=-3cm, align=center] {
        Compensation equation (\ref{eq:dots})\\ 
        Backward simulation to yield $\boldsymbol{s}(t)$ 
        };
        \node(qtut) [rec=gray, below of=WbQbM, yshift=-0.2cm, align=center]{
        ~\textcolor{black}{Integrate linear reduced dynamics (\ref{eq:state-qt-dynamics}) for $\boldsymbol{{q}}(t)$}\\
        and compute feedback control input $\boldsymbol{u}(t)$ (\ref{eq:u-closed-loop-q})
        };
        \node(zt) [rec=gray, below of=qtut, yshift=1.5cm, align=center] {
        Predicted response $\boldsymbol{z}(t)=\boldsymbol{W}(\boldsymbol{p}(t))+\varepsilon\boldsymbol{\hat{V}q}(t)$ (\ref{eq:decomp-zt-qt}) 
        };
        
    \draw[arrow](SSM.south) --node[left]{}($(SSM.south)+(0,-0.3cm)$) --($(IC.north)+(0,0.4cm)$) --(IC.north);
    \draw[arrow](uv.south) --($(uv.south)+(0,-0.3cm)$) --($(IC.north)+(0,0.4cm)$) --(IC.north);
    \draw[arrow](IC.south) --(calP.north);
    \draw[arrow](IC.east)  --($(Ricti.north)+(0,3.78cm)$) --(Ricti.north);
    \draw[arrow](calP.south) --(WbQbM.north);
    \draw[arrow](WbQbM.west) --($(st.north)+(0,0.85cm)$) --(st.north);
    \draw[arrow](Ricti.west) --(st.east);

    \draw[dashed, black] ($(Ricti.south)+(0cm,-0.1cm)$) -- ++(3cm,0cm);
    \draw[dashed, black] ($(Ricti.south)+(0cm,-0.1cm)$) -- ++(-3cm,0cm);
    \draw[dashed, black] ($(st.south)+(0cm,-0.1cm)$) -- ++(-3cm,0cm);
    \draw[dashed, black] ($(st.south)+(0cm,-0.1cm)$) -- ++(3cm,0cm);

    \draw[arrow](Ricti.south) --($(Ricti.south)+(0,-0.5)$) --($(qtut.north)+(0,0.3cm)$) --($(qtut.north)$);
    \draw[arrow](st.south) --($(st.south)+(0,-0.5)$) --($(qtut.north)+(0,0.3cm)$) --($(qtut.north)$);
    \draw[arrow](qtut.south) --(zt.north);

    \end{tikzpicture}
\caption{A flowchart depicting the implementation of our control design framework via SSM-based model reduction. Here, the online computation is minimal because only the simulation of the linear dynamics for the reduced coordinates $\boldsymbol{q}$ is involved in the online stage. In some of our examples, the dimension of the vector $\boldsymbol{q}$ is just two. The offline computation cost is also small because the dimension of vector $\boldsymbol{p}$ often equals to two or four.}
\label{fig:flowchart}
\end{figure*}

We conclude this subsection with a flowchart that summarizes the key ingredients of the above procedure of our closed-loop solution to the ELQP (\ref{P4}). The flowchart is presented in Fig.~\ref{fig:flowchart}. We note that the online computation is minimal because only the simulation of a linear ODE for the reduced linear coordinates $\boldsymbol{q}$ is involved, paving the road for the real-time active vibration control of flexible structures. Indeed, as we will see in later examples, the dimension of the reduced linear coordinates $\boldsymbol{q}$ is often equal to two, which is essentially a minimal system with a single degree of freedom. We also note that the offline computational cost is small because the SSMs are often two or four-dimensional and can be efficiently computed via SSMTool, as illustrated in~\cite{jain2022compute,part-i}.

\subsection{Receding horizon control}
\label{sec:receding}
The effectiveness of our SSM-based control design can deteriorate if the time horizon $[t_0,t_f]$ is long. Indeed, we observe in some examples considerable discrepancies between the predicted response and the actual response with the control law applied when the horizon is very long. To resolve the discrepancies, we extend our control design to receding horizon control~\cite{Mattingley2011RecedingHC}. Specifically, we divide the whole time horizon into several subintervals and apply the SSM-based control design shown in Fig.~\ref{fig:flowchart} to each subinterval. At the intersection point of two successive subintervals, the true final state of the former subinterval is taken as the initial state of the later subinterval to feed the correction at the intersection point to the next subinterval. A remaining question is how to select the length of each subinterval. We observe that the vibration of an uncontrolled system often has a characteristic time scale, which could be the lowest natural frequency or the period of self-excited oscillation. Then, we can take a half-period or longer window as the length of subinterval to conduct the receding horizon control.

\section{Numerical examples}
\label{sec:numerialExamples}

We consider a suite of examples with increasing complexity to illustrate the effectiveness of the proposed control design via the SSM-based ROMs, including a finite element model of an aircraft wing with more than 13,000 degrees of freedom. In each example, we will validate the effectiveness of our control designs by comparing ROM-based prediction and the corresponding response of the full system. As we will see, our control design achieves desirable accuracy and remarkable efficiency, paving the way for real-time control of nonlinear vibration of high-dimensional flexible structures. We perform computations for the aircraft wing on a workstation with Intel Xeon Platinum 8160 CPU @2.10GHz and 1024 GB RAM. The rest of the computations are performed on a PC with Intel Core i9-12900K CPU @3.20 GHz and 64 GB RAM.

\subsection{An oscillator chain}

\begin{figure}[!ht]
    \centering
    \includegraphics[width=1\linewidth]{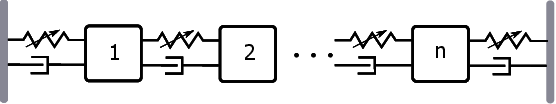}
    \caption{The schematic of an $n$-mass oscillator chain connected via cubic nonlinear springs. In this example, we set $n=10$.}
    \label{fig:oscillator chain}
\end{figure}

As our first example, we consider a chain of 10 oscillators with cubic nonlinearity shown in Fig.~\ref{fig:oscillator chain}. The oscillator chain can be interpreted as a discrete model of some flexible structures via lumped masses. The equations of the oscillator chain are in the form~\eqref{eq:eom-second-full} with $\boldsymbol{M}=m\boldsymbol{I}_{10}, \boldsymbol{K}=k\boldsymbol{T}, \boldsymbol{C}_\mathrm{d}=c\boldsymbol{T}$, and
\begin{align}
\label{eq:osc-f-E}
    \boldsymbol{f}(\boldsymbol{x},\dot{\boldsymbol{x}})&=\kappa
    \begin{bmatrix}
        x_1^3-(x_2-x_1)^3 \\
        (x_2-x_1)^3-(x_3-x_2)^3 \\
        \vdots \\
        (x_{10}-x_{9})^3-x_{10}^3 
    \end{bmatrix},\nonumber \\
\mathbf{E}(\Omega t)&=\boldsymbol{D}
\begin{pmatrix}
    \text{sin}(0.1\sqrt{2} t) \\ \text{cos}(0.1\sqrt{3} t)
\end{pmatrix}.
\end{align}
Here,
\begin{align}
    \boldsymbol{T}&=
    \begin{bmatrix}
        2  & -1 \\
        -1 &  2  & -1 \\
           & \ddots & \ddots & \ddots \\
        &  & -1 & 2 & -1 \\
        &  &  & -1 & 2
    \end{bmatrix}\in\mathbb{R}^{10\times 10},\nonumber \\
    \mathbf{D}&=
    \begin{bmatrix}
        1 & 0 & 0 & 0 & 0 & 0 & 0 & 0 & 0 & 0  \\
        0 & 0 & 0 & 0 & 1 & 0 & 0 & 0 & 0 & 0  
    \end{bmatrix}^T,
\end{align}
and $\boldsymbol{u}(t)\in{\mathbb{R}^2}$ is a two-dimensional control vector. In particular, we place the two actuators at the first and fifth oscillators. System parameters are chosen as $m=1,k=1,c=0.1,\kappa=0.5$ and $\epsilon=0.001$ in the following computations.

With the above parameter values, the pair of complex conjugate eigenvalues with the largest real parts is obtained as $\lambda_{1,2}=-0.0041\pm 0.2846\mathrm{i}$. The associated pair of modes is selected to span the master spectral subspace for SSM-based model reduction. We compute the autonomous SSM and its associated reduced dynamics at $\mathcal{O}(3)$ using SSMTool~\cite{ssmtool21}. The expansion order is determined via the convergence of backbone curves, as detailed in Remark~\ref{rk1}.

\begin{figure}[!ht]
    \centering
    \includegraphics[width=0.5\textwidth]{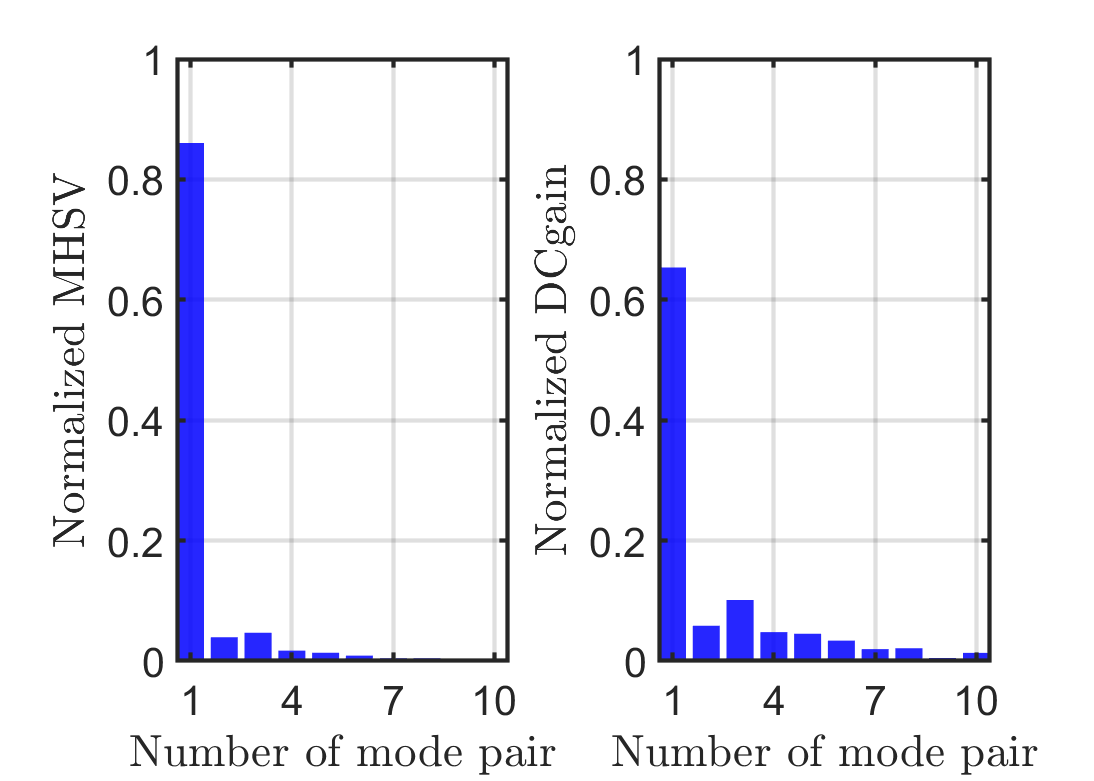}
    \caption{Normalized DCgain and MHSV for a system of oscillator chain. Here and throughout this paper, the `number of mode pair' in horizontal axis denote the index of mode pair instead of the total number of mode pairs.}
    \label{fig:DC_OscillatorChain}
\end{figure}

Next, we select the reduction basis $\hat{\boldsymbol{V}}$ for the linear dynamics with control. Here the displacements of the first and fifth oscillators are taken as observable. Following Sect.~\ref{sec:selection-basis}, we compute DCgains and MHSVs for the first ten pairs of eigenmodes, and the obtained results are plotted in Fig.~\ref{fig:DC_OscillatorChain}. We find that the summation of normalized DCgains and MHSVs for the first five pairs of modes is equal to 90.7\% and 97.8\%, respectively, suggesting that the reduction on the first five pairs of modes is of sufficient accuracy. Thus, we use these five pairs of modes to construct $\hat{\boldsymbol{V}}$ for the purpose of linear reduction.

\begin{figure}[!ht]
    \centering
    \includegraphics[width=0.475\textwidth]{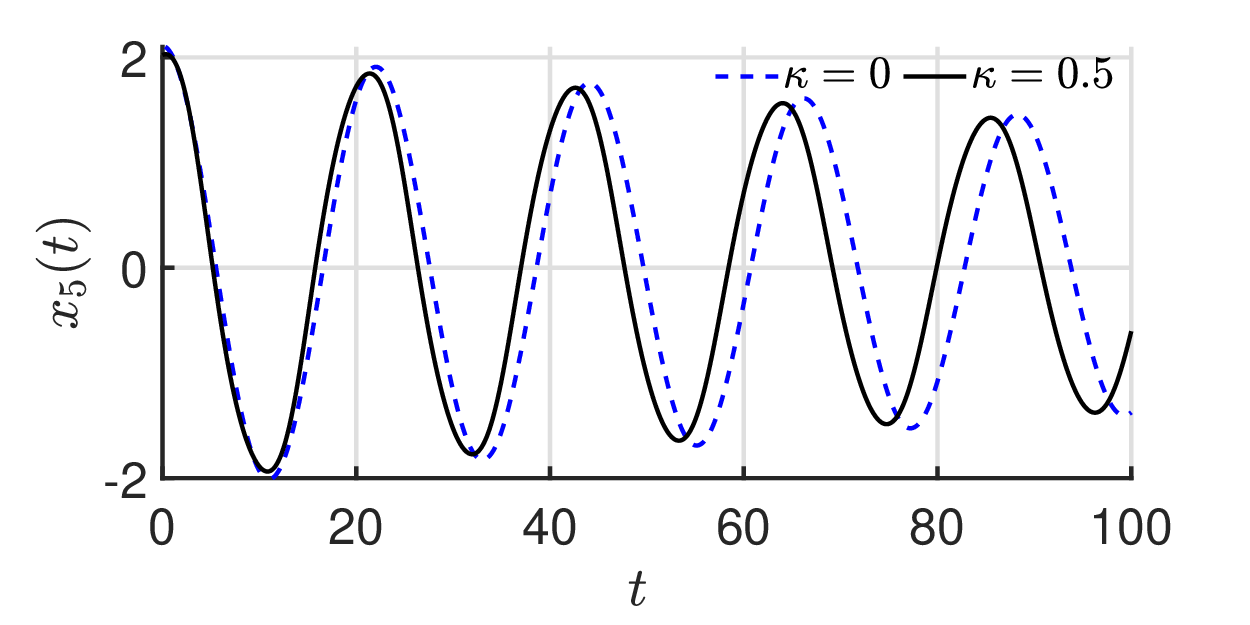}
    \caption{Free oscillation of full model for the oscillator chain. The black solid line gives the nonlinear free oscillation and the blue dashed line represents the linear free oscillation.}
    \label{fig:OscillatoChain_NL}
\end{figure}

We set the initial value of reduced coordinates $\boldsymbol{p}_0=[2.5,2.5]^\text{T}$, take $\boldsymbol{z}_0=\boldsymbol{W}(\boldsymbol{p}_0)$ as the initial conditions of the full system. In particular, the corresponding initial displacement of the fifth mass is 2.0217. This initial displacement is large enough such that the cubic nonlinearity in~\eqref{eq:osc-f-E} affects the vibration, as suggested by Fig.~\ref{fig:OscillatoChain_NL}, where the free vibration with the cubic nonlinearity ($\kappa=0.5$) and without it ($\kappa=0$) is shown.

\begin{figure}[!ht]
    \centering
    \includegraphics[width=0.5\textwidth]{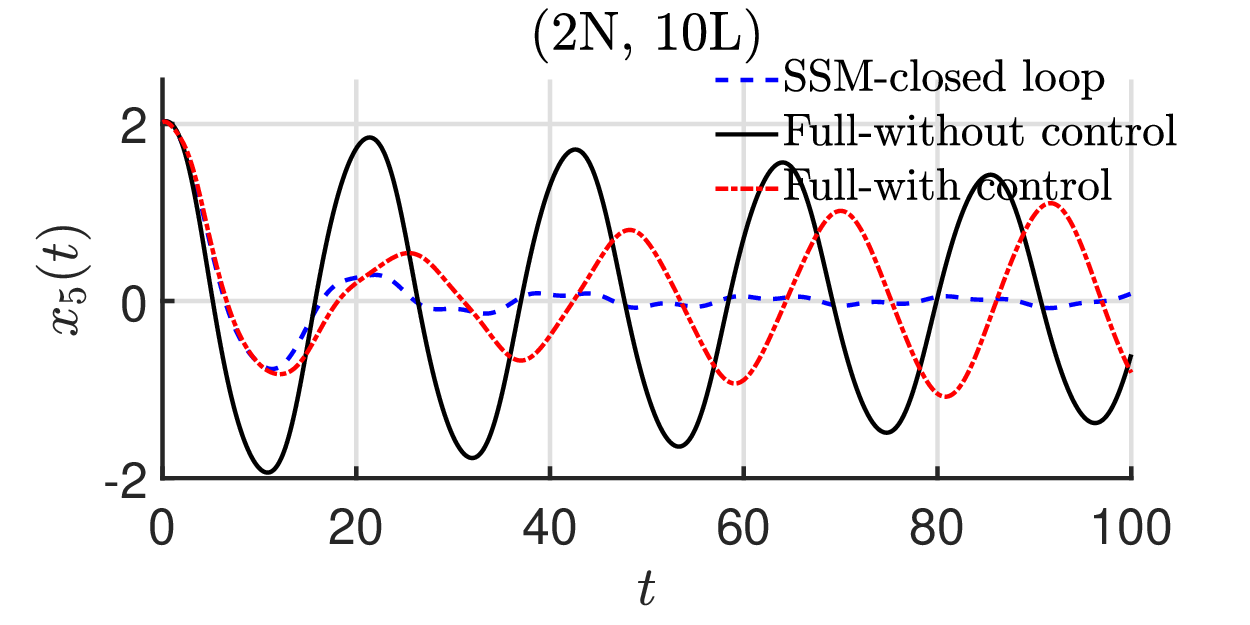}
    \includegraphics[width=0.5\textwidth]{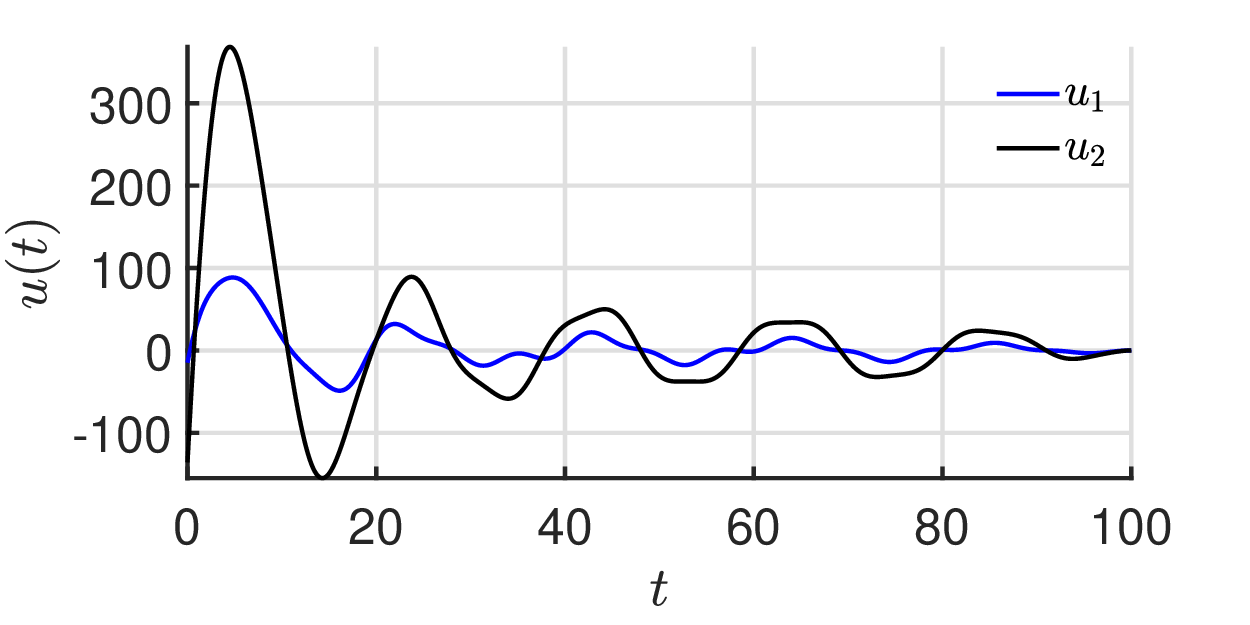}
    \caption{Vibration control of an oscillator chain with two actuators. In the upper panel, the response of the fifth oscillator is given, where the black solid line denotes the response of the full system without control, the blue dashed line represents the controlled response predicted via SSM-based model reduction, and the red dashed-dotted line gives the response of the full system with the control applied. In the lower panel, the control signals of the two actuators obtained using our SSM-based model reduction are given. Here and throughout this paper, $(m_1\mathrm{N},m_2\mathrm{L})$ stands for a $m_1$-dimensional nonlinear reduction and a $m_2$-dimensional linear reduction, namely, $\boldsymbol{p}\in\mathbb{C}^{m_1}$ and $\boldsymbol{q}\in\mathbb{C}^{m_2}$. In this example, we have $m_1=2$ and $m_2=10$.}
    \label{fig:OscillatoChain_LQR}
\end{figure}

Under these initial conditions and the external forcing, the transient response of the fifth oscillator for $t\in[0,100]$ is plotted in a black solid line in the upper panel of Fig.~\ref{fig:OscillatoChain_LQR}. We apply our control algorithm to suppress the large amplitude oscillation of the chain. The weight matrices of LQR are given as $\mathbf{Q}=\text{diag}(10^5\mathbf{I}_{10},0,\cdots,0), \mathbf{\hat{R}}=\text{diag}(0.05,0.05)$ and $\mathbf{\hat{M}}=\mathbf{0}$. We solve the ELQP~\eqref{P4} following the procedure shown in Fig.~\ref{fig:flowchart}. The control input produced via the ROM-based control algorithm is shown in the lower panel of Fig.~\ref{fig:OscillatoChain_LQR}, and the corresponding controlled response for $x_1(t)$ is plotted in a blue dashed line in the upper panel of Fig~\ref{fig:OscillatoChain_LQR}. The vibration is suppressed quickly and then restricted to small amplitude oscillations with the control applied.

To validate the effectiveness of the control design, we apply the control signals predicted via the ROM-based algorithm to the full system and compare the response of the full system with that of the ROM-based prediction. Specifically, we simulate the full system with the control input shown in the lower panel of Fig.~\ref{fig:OscillatoChain_LQR}. The simulated response of the full system is plotted in red dashed-dotted line in the upper panel of Fig.~\ref{fig:OscillatoChain_LQR}. We observe that the simulated response of the full system matches well with that of the ROM-based prediction for $t\in[0,20]$, but starts to deteriorate for $t\geq20$. This deterioration results from the accumulation of errors as the time horizon increases.

\begin{figure}[!ht]
    \centering
    \includegraphics[width=0.5\textwidth]{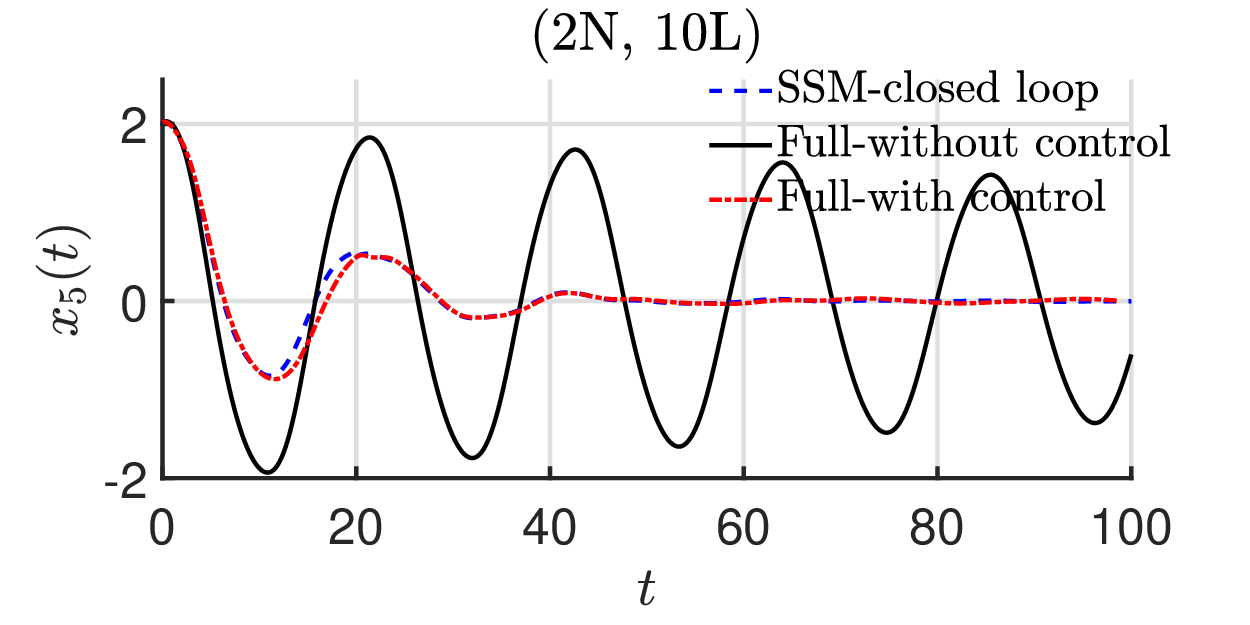}
    \includegraphics[width=0.5\textwidth]{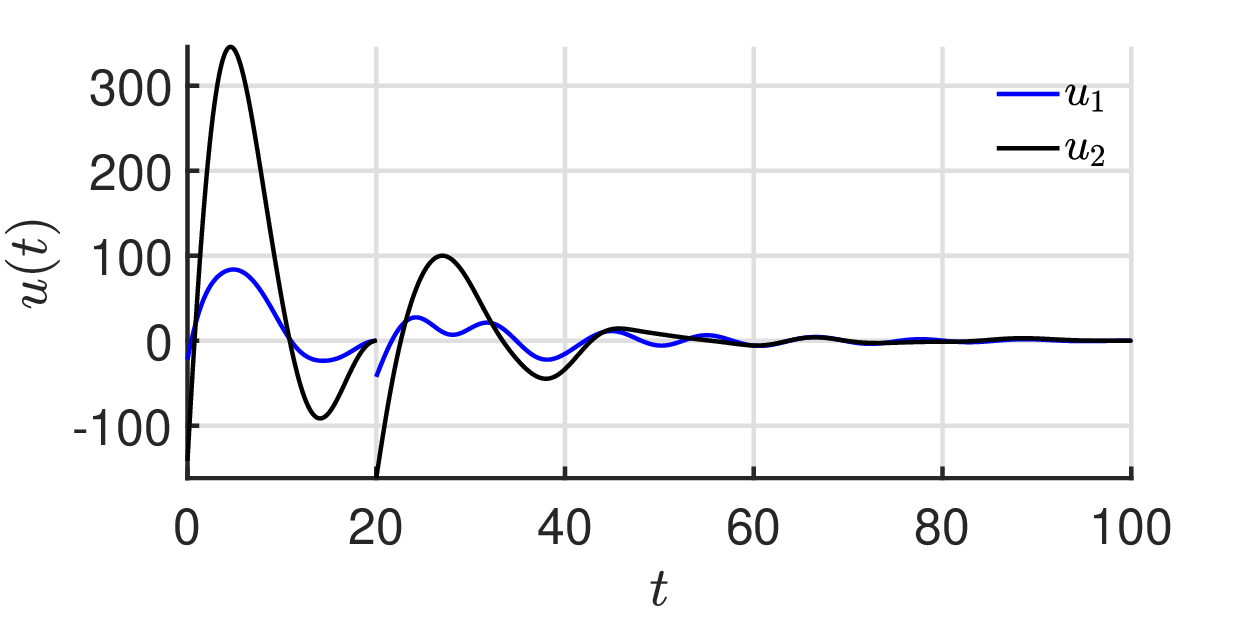}
    \caption{Vibration control of the oscillator chain with nonlinear oscillation. Unlike the case of Fig.~\ref{fig:OscillatoChain_LQR}, here, receding horizon control is applied. In particular, the whole time horizon is divided into two subintervals with the intersection point at $t=20$.}
    \label{fig:OscillatoChain_LQR_feedback}
\end{figure}

To address the deterioration, we apply the receding horizon control established in Sect.~\ref{sec:receding} to this example. In particular, we divide the horizon $[0,100]$ into two subintervals with an intersection point at $t=20$. The control inputs for the two subsequent subintervals $[0,20]$ and $[20,100]$ are shown in the lower panel of Fig.~\ref{fig:OscillatoChain_LQR_feedback}. We observe that the control signals are discontinuous at $t=20$, which is the intersection point of the two subintervals. This discontinuity highlights the correction by the feedback at the intersection point. Importantly, we see that the response of the full system matches excellently with the ROM-based prediction at the second subinterval. This illustrates the effectiveness of the receding horizon control.

\subsection{A von K{\'a}rm{\'a}n beam}

As our second example, we consider the vibration control of a clamped-pinned beam with von K{\'a}rm{\'a}n type of geometric nonlinearity. The beam's length, width, and height are 2700 mm, 10 mm, and 10 mm, respectively. Material properties are given by a density $\rho=1780\times10\text{kg/m}^3$, Young's modulus $E=45\text{GPa}$, and Poisson's ratio $\nu=0.3$. The beam is divided into 8 elements, resulting in a system with 22 degrees of freedom. We adopt Rayleigh damping $\boldsymbol{C}_\mathrm{d}=\alpha{\boldsymbol{M}}+\beta{\boldsymbol{K}}$ with $\alpha=0$ and $\beta={1/450}~\mathrm{s}^{-1}$.

The first two pairs of complex conjugate eigenvalues of the finite element (FE) model of the beam are obtained as $\lambda_{1,2}=-1.0472\pm30.6815\mathrm{i}$ and $\lambda_{3,4}=-11.0053\pm 98.9126\mathrm{i}$. Since $\lambda_3\approx3\lambda_1$, this system admits a near 1:3 internal resonance between the first two pairs of modes~\cite{Nayfeh1974NonlinearAO}. Thus, we take their associated eigenvectors to span a four-dimensional master spectral subspace for SSM-based model reduction. We compute the autonomous SSM and its reduced dynamics in $\mathcal{O}({5})$ expansion. This expansion order is again determined via the procedure detailed in Remark~\ref{rk1}.

\begin{figure}[!ht]
    \centering
    \includegraphics[width=0.5\textwidth]{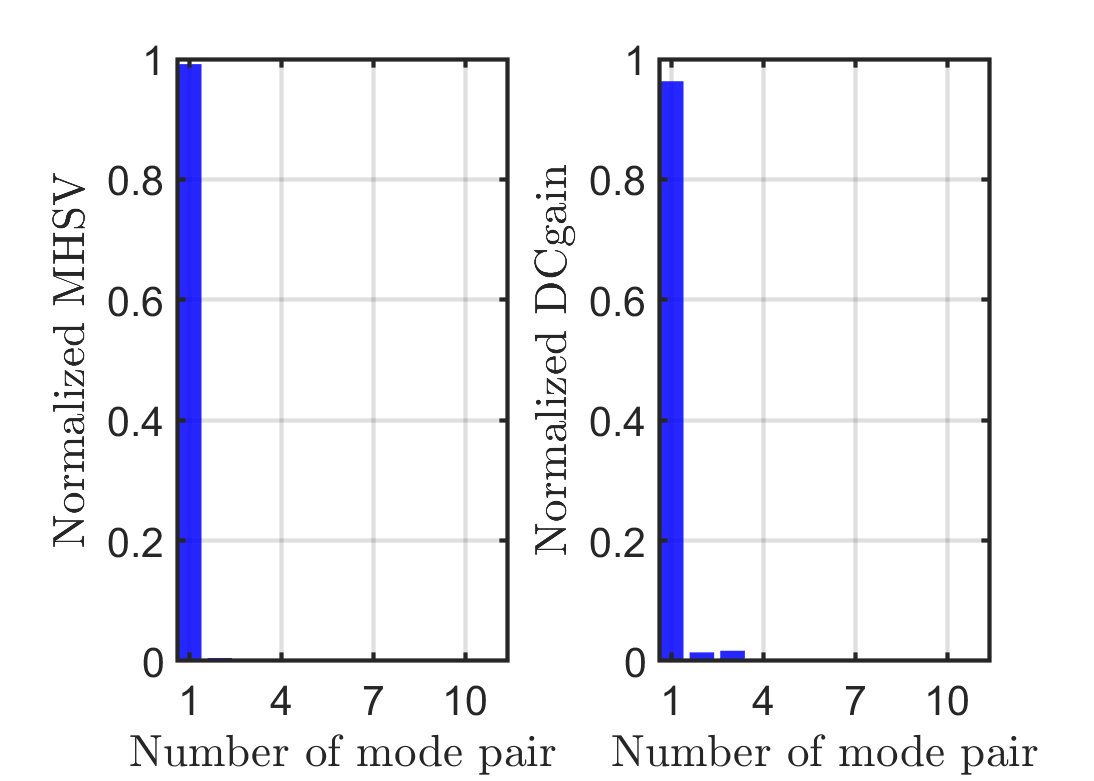}
    \caption{Normalized DCgain and MHSV for a von K{\'a}rm{\'a}n beam.}
    \label{fig:VonKarmanBeamDCandHSV}
\end{figure}

Next, we select the reduction basis $\hat{\boldsymbol{V}}$ for the linear dynamics with control. Following Sect.~\ref{sec:selection-basis}, we compute DCgains and MHSVs for the first ten pairs of eigenmodes, and the obtained results are plotted in Fig.~\ref{fig:VonKarmanBeamDCandHSV}. Here, an actuator is placed at the midpoint of the beam, and the deflection at the midpoint is taken as observable. The DCgain and MHSV for the first pair of modes is equal to 96.3\% and 99.1\%. Therefore, using only the first pair of modes to perform the linear reduction is sufficient.

\begin{figure}[!ht]
    \centering
    \includegraphics[width=0.475\textwidth]{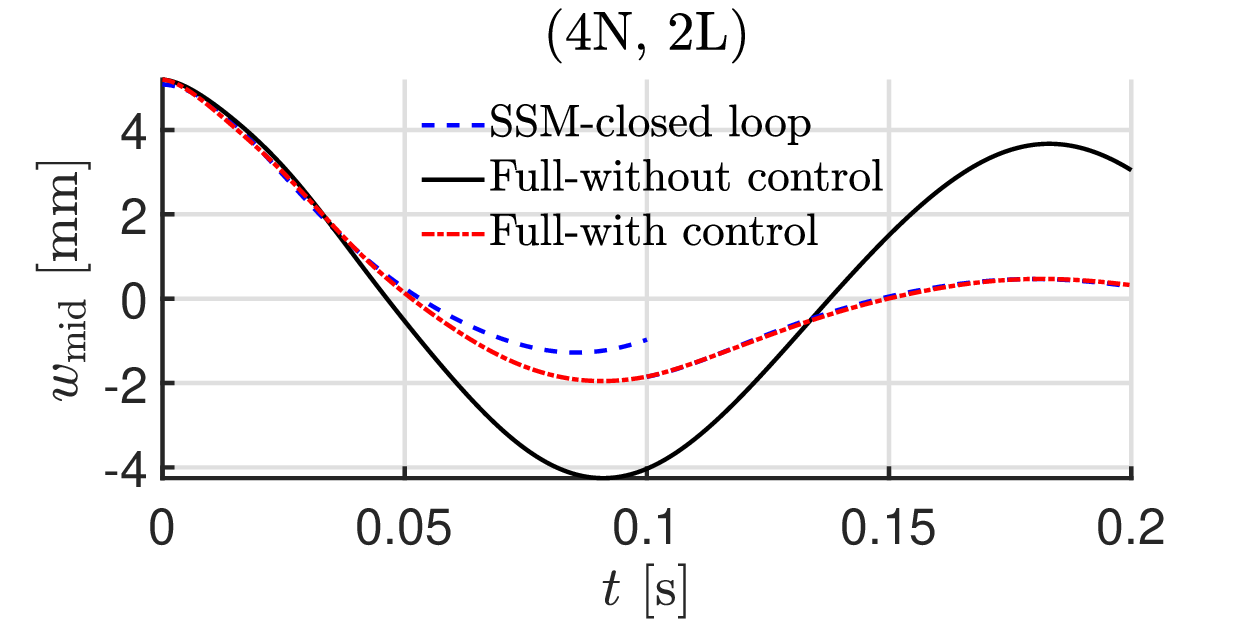}
    \includegraphics[width=0.475\textwidth]{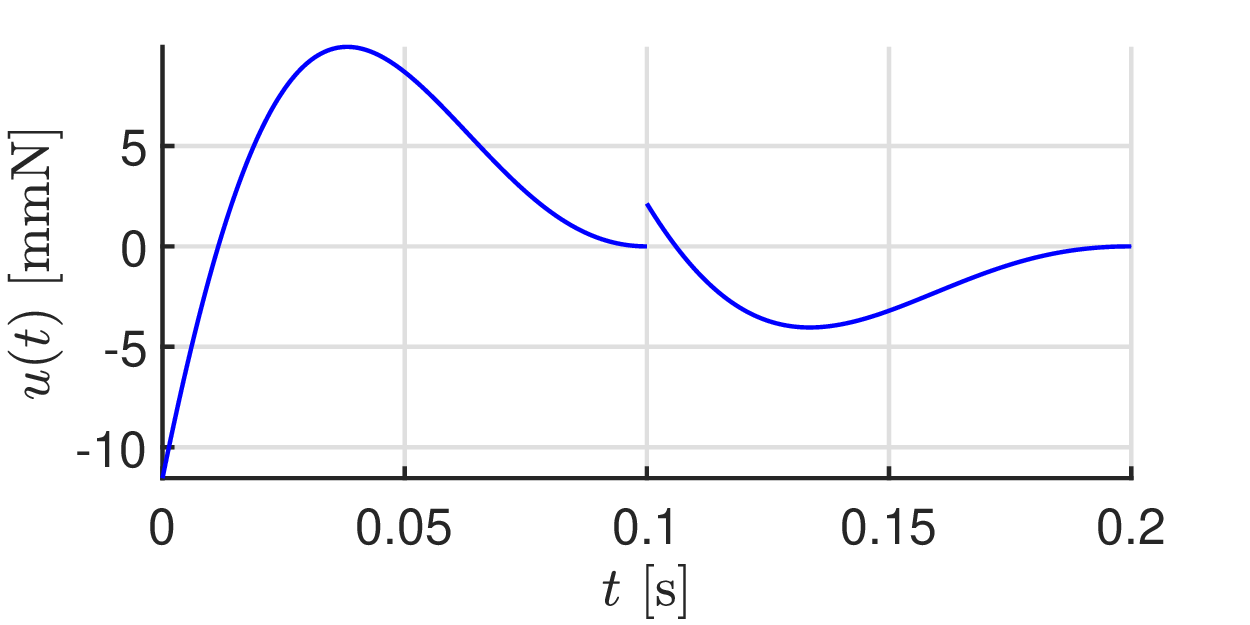}
    \caption{Vibration control of the von K{\'a}m{\'a}rn beam with an actuator placed at the midpoint of the beam. In the upper panel, the deflection at the midpoint of the beam is given, where the black solid line denotes the response of the full system without control, the blue dashed line represents the controlled response predicted via SSM-based model reduction, and the red dashed-dotted line gives the response of the full system with the control applied. In the lower panel, the control signal of the actuator obtained using our SSM-based model reduction is given. Here, receding horizon control is applied. In particular, the whole time horizon is divided into two subintervals with the intersection point at $t=0.1$.}
    \label{fig:VonKarmanBeamLQR2segs}
\end{figure}

We apply a transverse static loading of 1500 N at the midpoint of the beam, resulting in a deflection of 5 mm at the midpoint of the beam, which is half of the thickness of the beam. Therefore, the geometric nonlinearity will play an important role. We then release the static loading, and the beam undergoes a free vibration. Let $w_\mathrm{mid}$ be the deflection at the midpoint of the beam, and the free vibration regarding $w_\mathrm{mid}$ for $t\in[0,0.2]$ is plotted in black solid line in the upper panel of Fig.~\ref{fig:VonKarmanBeamLQR2segs}. Likewise, we apply our ROM-based control algorithm to suppress the vibration. With $\mathbf{Q}=\text{diag}(0.1\mathbf{I}_{22},\mathbf{0}_{22}),{\hat{R}}=0.01,\mathbf{\hat{M}}=\mathbf{0}$ and $\epsilon=50$, we apply the receding control algorithm with an intersection point at $t=0.1$ to suppress the large amplitude nonlinear vibration. The obtained control input is shown in the lower panel of Fig.~\ref{fig:VonKarmanBeamLQR2segs}, and the associated controlled response is presented in the upper panel (blue dashed line). We observe that the vibration is indeed effectively suppressed with the control applied. The discontinuity shown in Fig.~\ref{fig:VonKarmanBeamLQR2segs} again illustrates the correction of receding control.

Similar to the previous example, we apply the control signals to the full system to check the effectiveness of the ROM-based control design. The simulated trajectory of the full system is plotted in a red dashed-dotted line in the upper panel of Fig.~\ref{fig:VonKarmanBeamLQR2segs}. Although small discrepancies are observed in the first segment for $t\approx0.1$, the simulated response of the full system matches excellently with that of ROM-based prediction in the second segment. This again demonstrates the effectiveness of our SSM-based control design.

\subsection{A von K{\'a}rm{\'a}n shell}
\begin{figure}[!ht]
    \centering
    \includegraphics[width=1\linewidth]{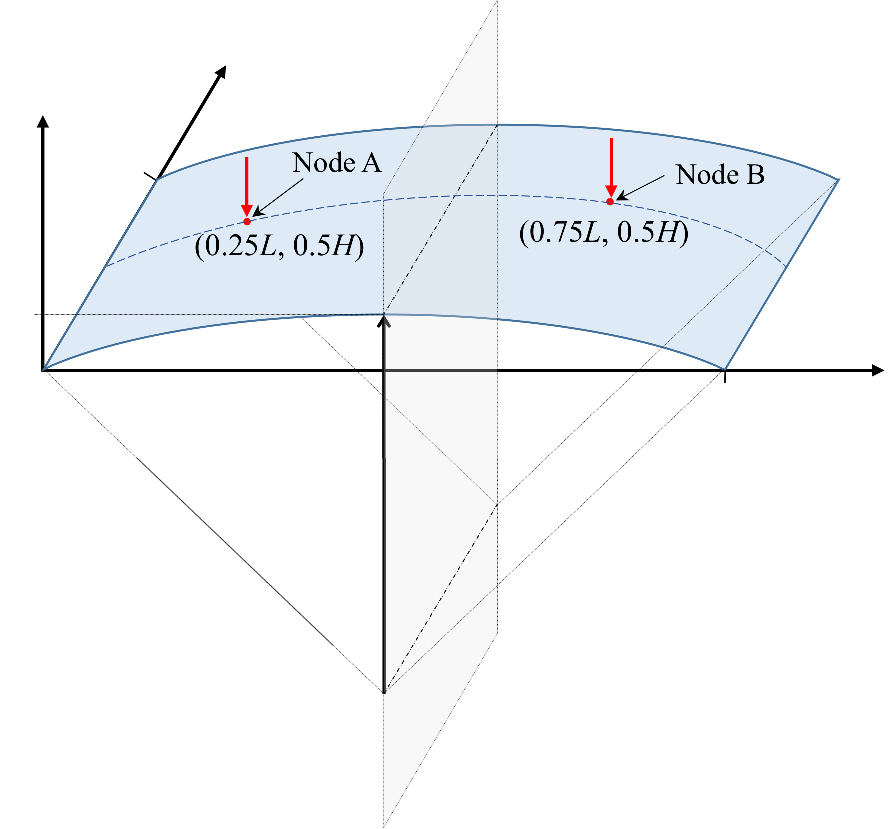}
    \caption{Schematic of a von K{\'a}rm{\'a}n shell structure}
    \label{fig:SchematicShell}
\end{figure}

Next, we study the vibration control of the finite element model of a geometry nonlinear von K{\'a}rm{\'a}n shell. The shell is supported at the two opposite edges aligned along the $y$-axis, as seen in Fig.~\ref{fig:SchematicShell}. Geometric parameters of the shell are given by the length $l=1\text{m}$, the breadth $b=2\text{m}$, the thickness $t=0.01\text{m}$, and the curvature radius is $c=0.1\text{m}$ (without internal resonance) or $c=0.041\text{m}$ (with internal resonance). Material properties are given as the density $\rho=2700\text{kg/m}^3$, Young's modulus $E=70\text{GPa}$, and the Poisson's ratio $\nu=0.33$.

The shell is divided into 400 flat, triangular shell elements, resulting in an FE model with 1320 degrees of freedom. We adopt Rayleigh damping and choose the damping coefficient $\alpha=0.4022$ and $\beta=8.6312\times{10}^{-6}~\mathrm{s}^{-1}$ such that the damping ratios of the first two modes are equal to 0.002. Here, we consider two cases: one without internal resonance and one with a 1:2 internal resonance between the first two bending modes of the shell to illustrate the wide applicability of our control design framework. As we will see, the dimensions of both the nonlinear and linear reduction will increase when the internal resonance is considered. This increase is a result of modal interactions in internally resonant structures.

\begin{figure}[!ht]
    \centering
    \includegraphics[width=0.375\textwidth]{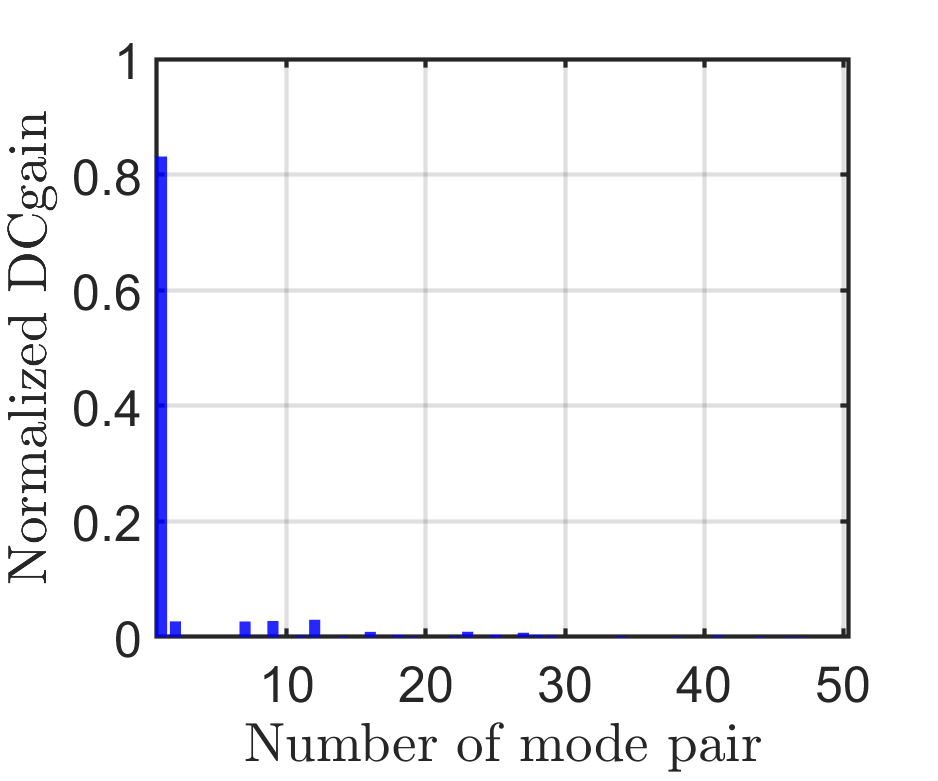}
    \includegraphics[width=0.375\textwidth]{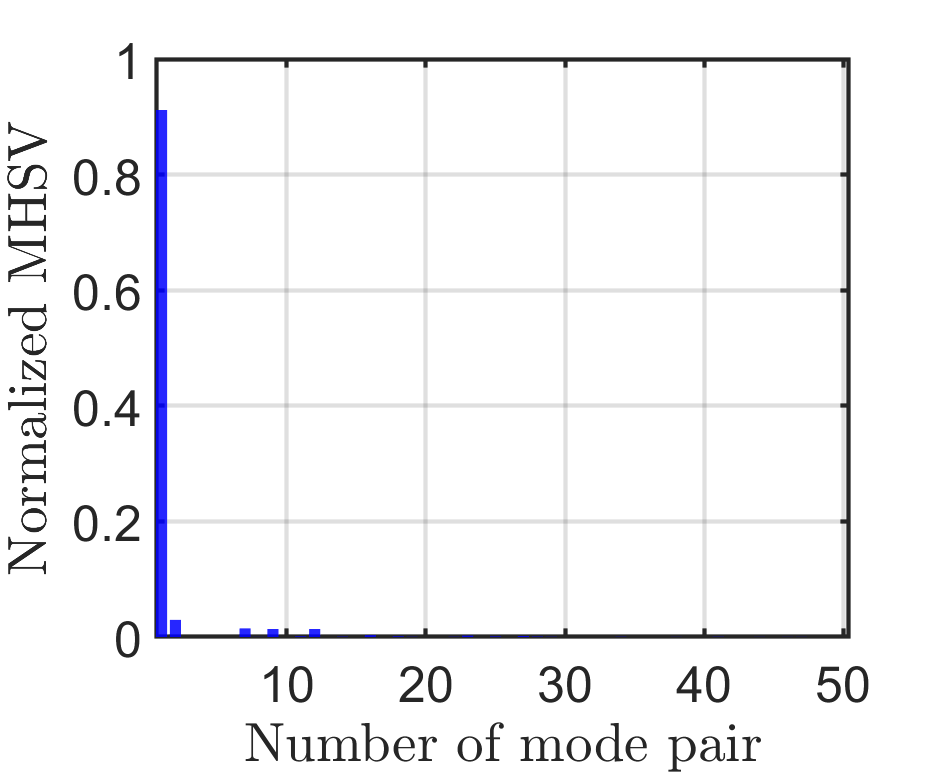}
    \caption{Normalized DCgain and MHSV for a von K{\'a}rm{\'a}n shell without internal resonance.}
    \label{fig:VonKarmanShellDCandHSV}
\end{figure}

We start with the case without internal resonance. The first pair of complex modes with associated eigenvalues $\lambda_{1,2}=-0.29\pm 147.45\mathrm{i}$ is chosen to span a two-dimensional master spectral subspace. We compute the autonomous SSM and its reduced dynamics at~$\mathcal{O}(5)$ expansion. We again use DCgain and MHSV to choose the reduction basis $\hat{\boldsymbol{V}}$ for the linear dynamics with control. Here we place two actuators at two points $(x,y)=(0.25L,0.5H)$ and $(x,y)=(0.75L,0.5H)$ and apply control input along $z$-axis. The associated two degrees of freedom are further taken as observable. Fig.\ref{fig:VonKarmanShellDCandHSV} shows the DCgains and MHSVs of the first 50 modes of the shell. The DCgain and MHSV for the first pair of modes equal $83.12\%$ and $91.18\%$. Thus, we only keep the first pair of modes to approximate the full model with $1320$ degrees of freedom, resulting in a significant reduction in the system's complexity. 

\begin{figure}[!ht]
    \centering
    \includegraphics[width=0.5\textwidth]{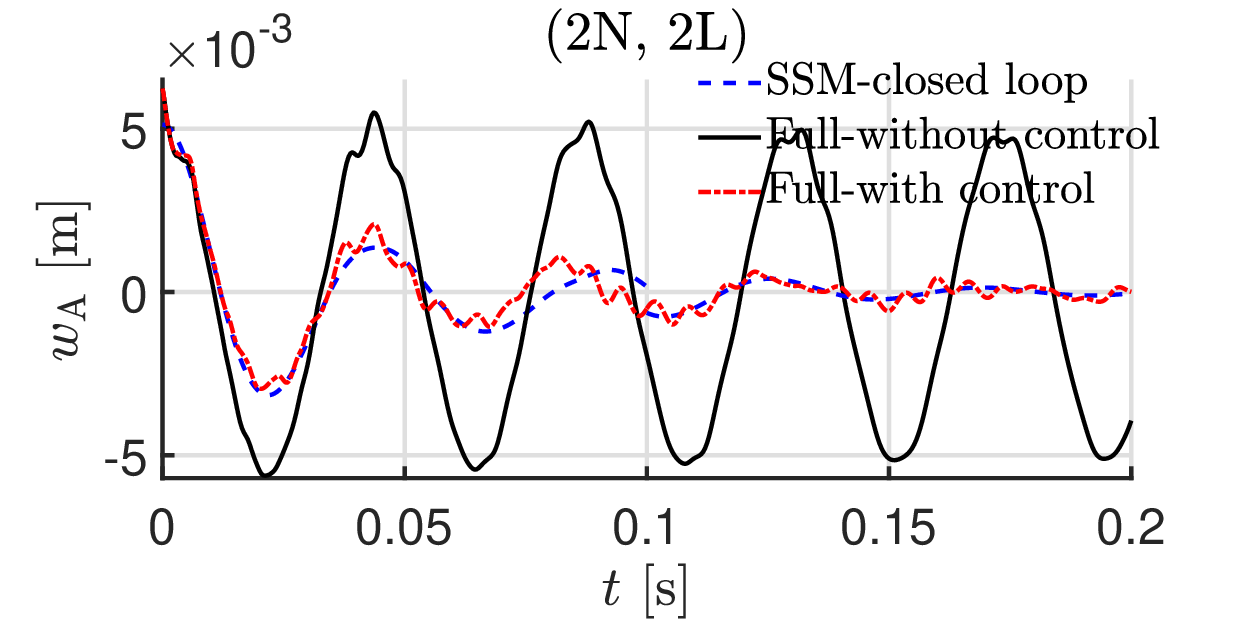}
    \includegraphics[width=0.5\textwidth]{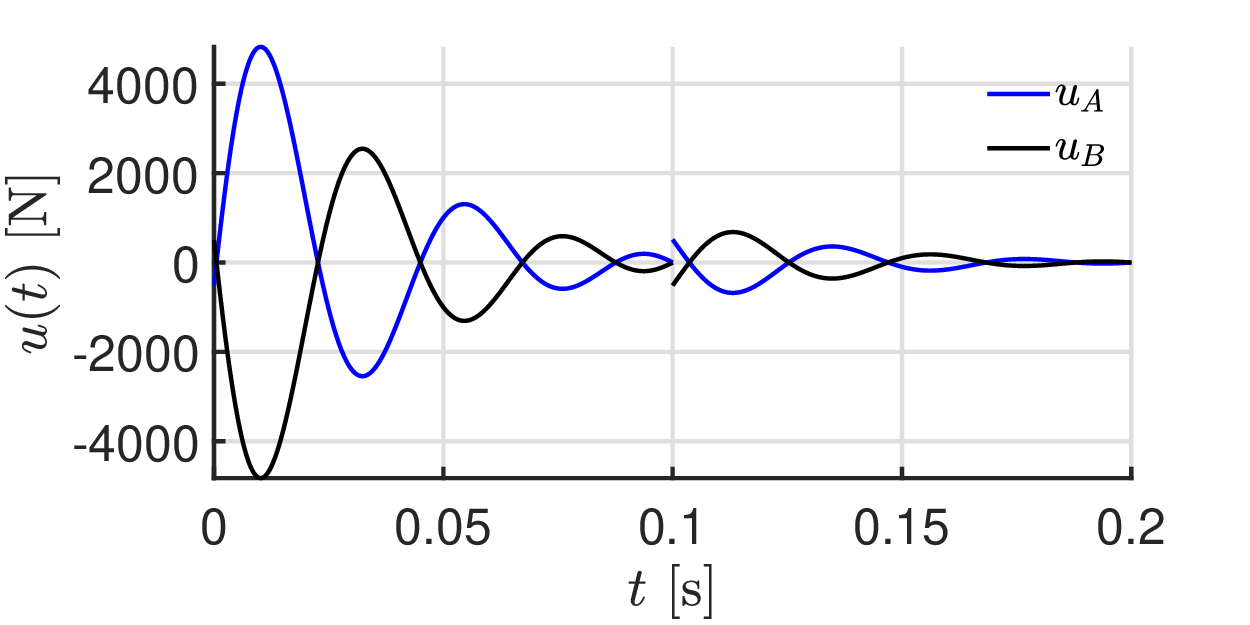}
    \caption{Vibration control of the von K{\'a}m{\'a}rn shell without internal resonance. Here, two actuators are placed on the two nodes shown in Fig.~\ref{fig:SchematicShell}, and the time horizon [0,0.2] is divided into two subintervals with an intersection point at $t=0.1$. In the upper panel, the deflection at node A is given, where the black solid line denotes the response of the full system without control, the blue dashed line represents the controlled response predicted via SSM-based model reduction, and the red dashed-dotted line gives the response of the full system with the control applied. In the lower panel, the control signals of the two actuators that were obtained using our SSM-based model reduction are given.}
    \label{fig:VonKarmanShellLQR2segs}
\end{figure}

Similar to the previous examples, we apply a transverse static loading 4000 N at node A of the shell, resulting in a deflection $6.2$ $\mathrm{mm}$ at node A, which is more than 60\% percent of the thickness of the shell. Therefore, the geometric nonlinearity of the shell will play a significant role. We release the static loading and then apply control to suppress this free vibration shown in the black solid line in Fig.~\ref{fig:VonKarmanShellLQR2segs}. With LQR weight matrices below
\begin{equation}
    \mathbf{Q}=10^4
    \begin{bmatrix}
        \mathbf{I}_n & \mathbf{0} \\
        \mathbf{0}  & \mathbf{I}_n 
    \end{bmatrix},
    \mathbf{\hat{R}}=
    \begin{bmatrix}
        0.1 & {0} \\
        {0}  & 0.1
    \end{bmatrix},
    \mathbf{\hat{M}}=\mathbf{0}.
\end{equation}
and $\epsilon=0.1$. We apply our control algorithm to the von K{\'a}rm{\'a}n shell with receding feedback at $t=0.1$s, and the obtained control result is shown in Fig.~\ref{fig:VonKarmanShellLQR2segs}. It suggests that the system is quickly stabilized after about two cycles of oscillation, and the reduced system approximates the full model well, as seen in the upper panel of the figure.

Now we move to the case with internal resonance. The first two pairs of eigenvalues are $\lambda_{1,2}=-0.30\pm149.22\mathrm{i}$ and $\lambda_{3,4}=-0.60\pm298.78\mathrm{i}$, suggesting there is a 1:2 internal resonance between the first two bending modes. This results in energy exchange between the two modes and can cause severe vibration and fatigue. Following SSM theory, we select the first two pairs of complex modes as the master subspace and compute the associated four-dimensional SSM and its reduced dynamics.

\begin{figure}[!ht]
    \centering
    \includegraphics[width=0.375\textwidth]{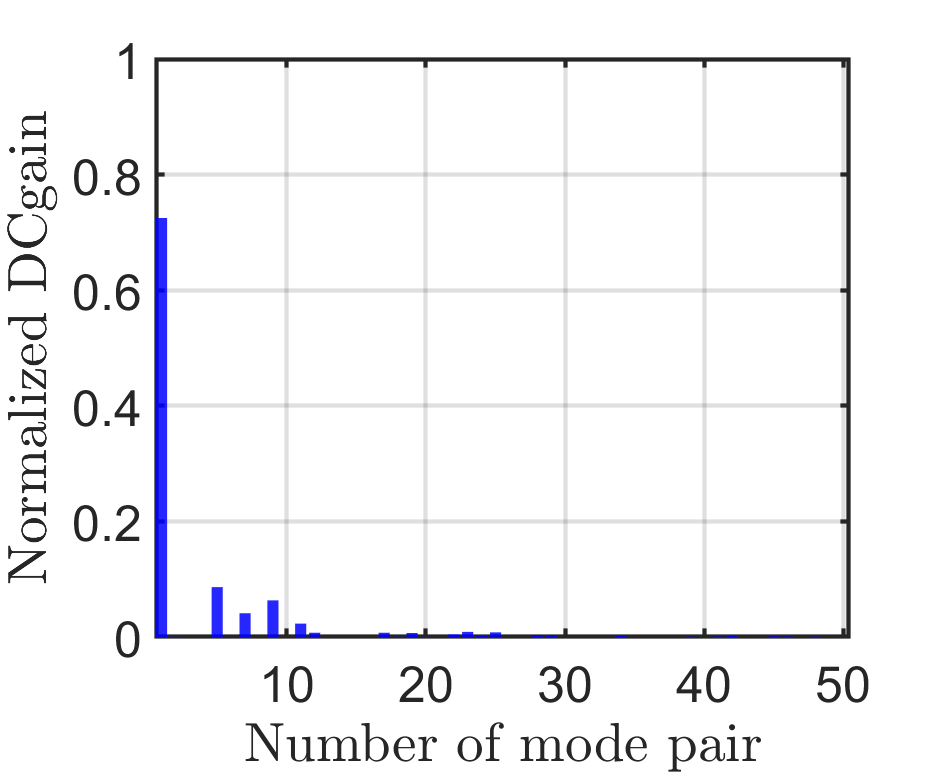}
    \includegraphics[width=0.375\textwidth]{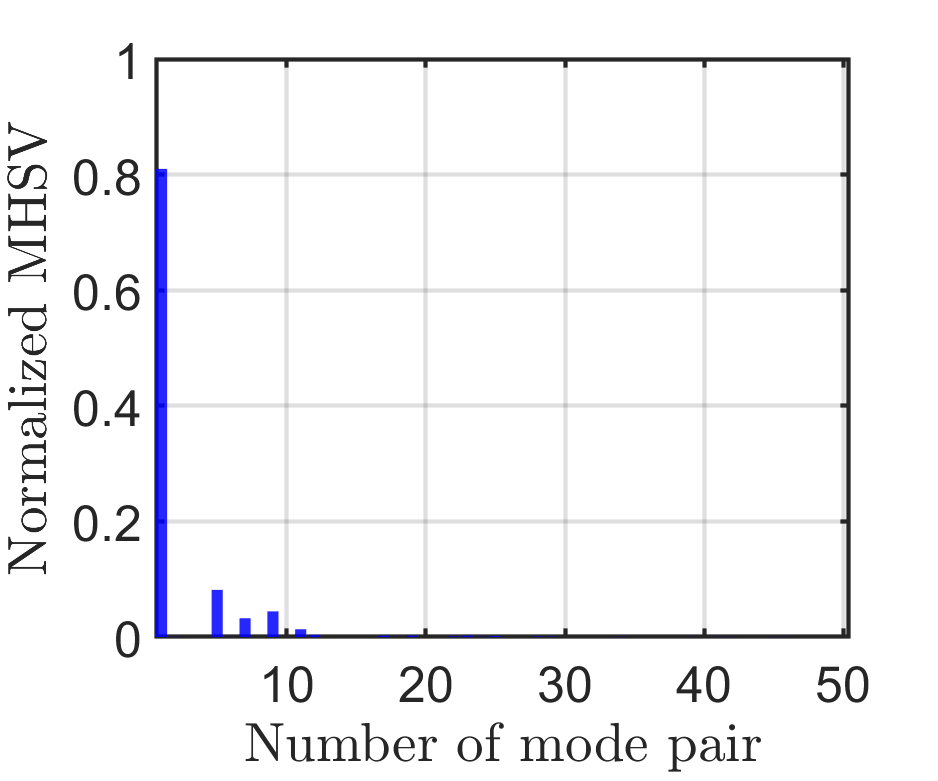}
    \caption{Normalized DCgain and MHSV for a von K{\'a}rm{\'a}n shell with 1:2 internal resonance.}
    \label{fig:IRs_VonKarmanShellDCandHSV}
\end{figure}

We again perform the computation of the DCgains and MSHVs for the first 50 pairs of eigenmodes and the obtained results are plotted in Fig.~\ref{fig:IRs_VonKarmanShellDCandHSV}. Unlike the case without internal resonance, where a linear reduction on a single pair of modes is sufficient, here we choose $1^\mathrm{st},5^\mathrm{th},7^\mathrm{th},9^\mathrm{th}$ pairs of eigenmodes to reduce the high-dimensional linear system because the summation of their normalized contributions is more than 97\% for DCgain and 91\% for MHSV. 

\begin{figure}[!ht]
    \centering
    \includegraphics[width=0.5\textwidth]{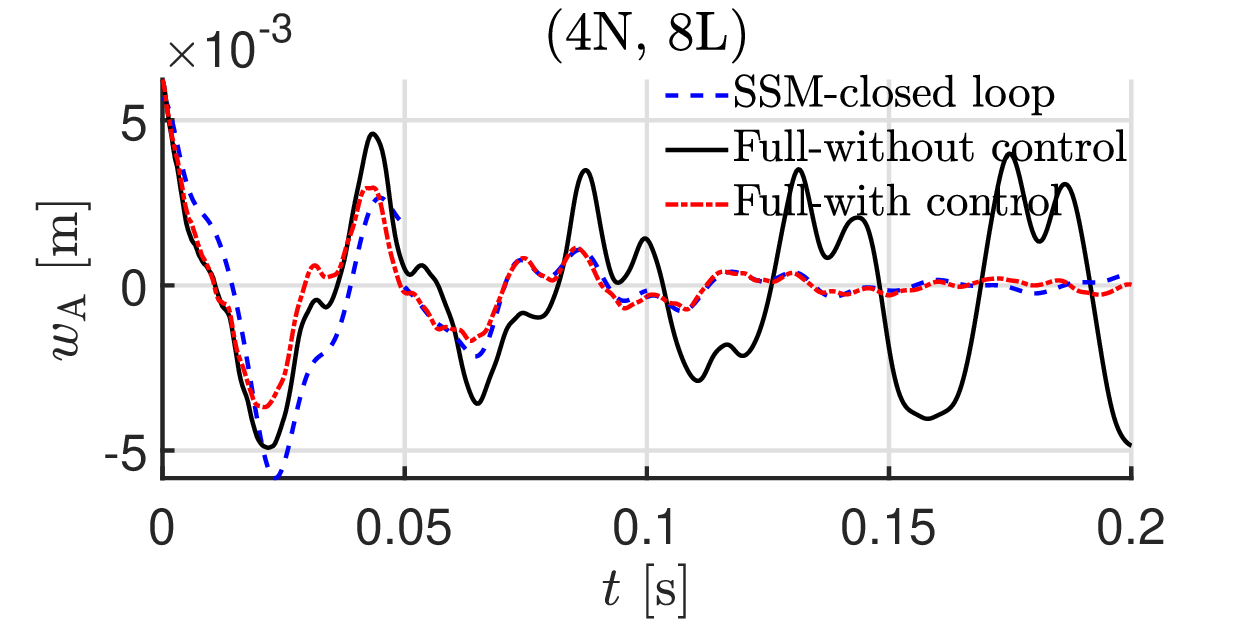}
    \includegraphics[width=0.5\textwidth]{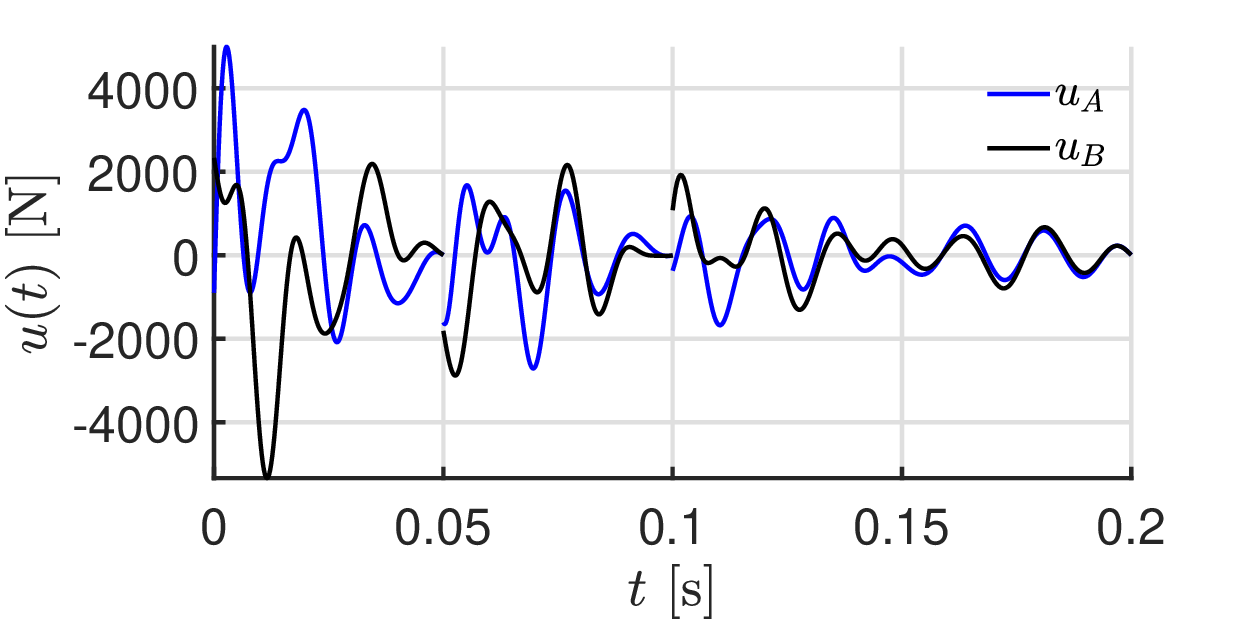}
    \caption{Vibration control of the von K{\'a}m{\'a}rn shell with a 1:2 internal resonance. Here, two actuators are placed on the two nodes shown in Fig.~\ref{fig:SchematicShell}, and the time horizon [0,0.2] is divided into three subintervals [0,0.05], [0.05,0.1] and [0.1,0.2]. In the upper panel, the deflection at node A is given, where the black solid line denotes the response of the full system without control, the blue dashed line represents the controlled response predicted via SSM-based model reduction, and the red dashed-dotted line gives the response of the full system with the control applied. In the lower panel, the control signals of the two actuators that were obtained using our SSM-based model reduction are given.}
    \label{fig:IRsVonKarmanShellLQR2segs}
\end{figure}

As seen in Fig.~\ref{fig:IRsVonKarmanShellLQR2segs}, the free vibration of the shell is of large amplitude. In particular, a growth regarding the response amplitude as time increases is observed for $t\approx0.15$. This growth further demonstrates the energy transfer due to modal interactions~\cite{liModelReductionConstrained2023}. We divide the time horizon $[0,0.2]$ into three subintervals and apply our control design framework to suppress the large amplitude nonlinear vibration. We observe from the upper panel of Fig.~\ref{fig:IRsVonKarmanShellLQR2segs} that the vibration is successfully suppressed and the controlled response predicted via the SSM-based model reduction match well with that of the full system, which illustrates the effectiveness of the control design. Comparing the lower panels of Fig.~\ref{fig:VonKarmanShellLQR2segs} and Fig.~\ref{fig:IRsVonKarmanShellLQR2segs}, we find that the control inputs for the case with internal resonance change more dramatically, to account for the contributions of higher-order modes.

\subsection{NACA wing}
As our fourth example, we consider the FE model of an aircraft wing with geometric nonlinearity, which was originally presented by Jain et al.~\cite{jainQuadraticManifoldModel2017}. The wing is clamped at one side, and the other side is free. Material parameters of the structure are given as Young's modulus $E=70$ Gpa, density $\rho=2700\text{kg/m}^3$, Poisson's ratio $\nu=0.33$. The thickness of the shell element is $t=0.0015$ m. More details on the geometric parameters can be found in~\cite{jainQuadraticManifoldModel2017}.

\begin{figure}[!ht]
    \centering
    \includegraphics[width=1\linewidth]{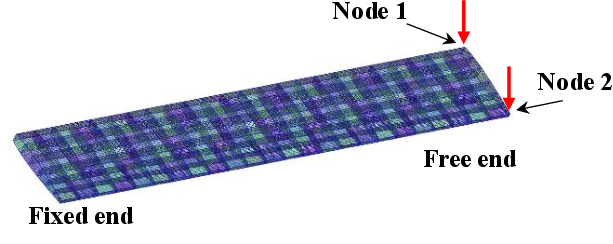}
    \caption{Schematic of FEM of an aircraft wing with 49968 flat triangular shell elements and 133920 degrees of freedom}
    \label{fig:NACAWingSchematic}
\end{figure}

The structure is meshed by 49968 flat triangular shell elements, resulting $n=133920$ degrees of freedom, as seen in Fig.\ref{fig:NACAWingSchematic}. Such a high dimensionality makes it impossible to control the system in real-time using the full model. We apply our SSM-based control design framework to suppress the nonlinear vibration of the wing. In this example, we again adopt Rayleigh damping such that a damping ratio of 0.4\% along the first two vibration modes is ensured. With linear spectral analysis, we choose the first pair of modes with eigenvalues $\lambda_{1,2}=-0.0587\pm29.3428\mathrm{i}$ to span the master spectral subspace for computing autonomous SSM.

\begin{figure}[!ht]
    \centering
    \includegraphics[width=0.8\linewidth]{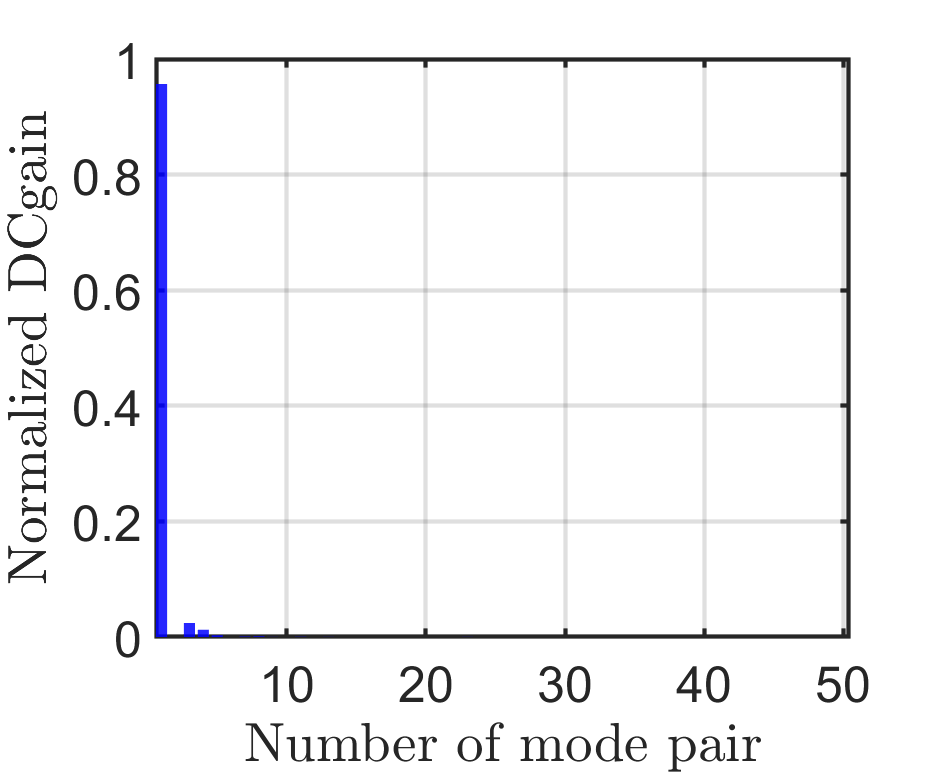}
    \includegraphics[width=0.8\linewidth]{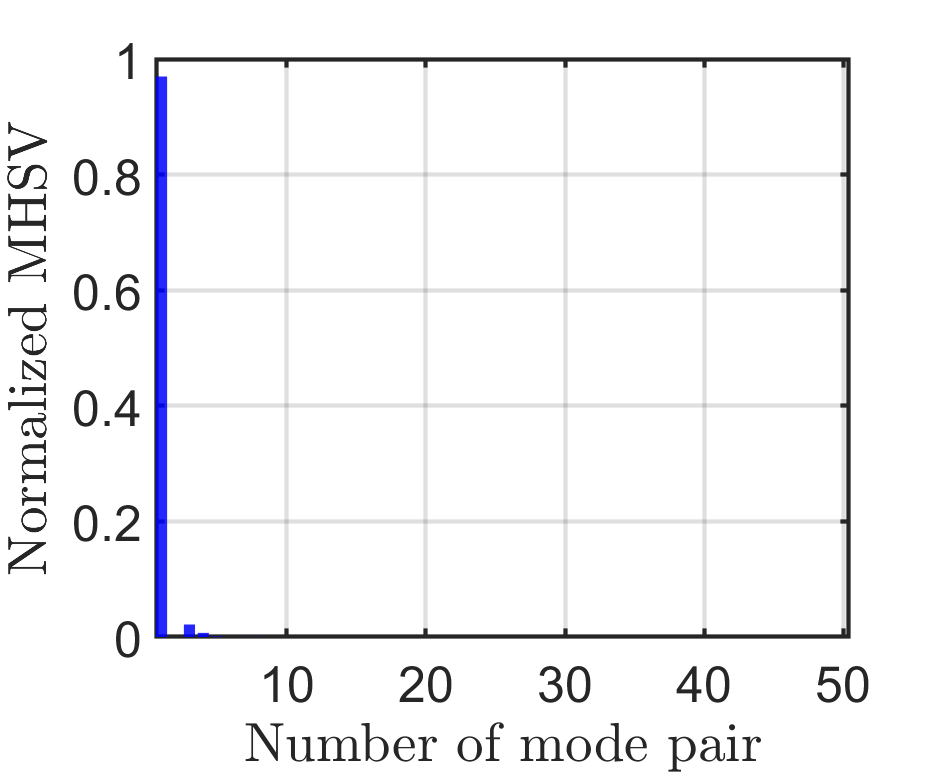}
    \caption{Normalized DCgain and MHSV for the first 50 pairs of modes of a NACA wing.}
    \label{fig:NACAWingHSVDCgain}
\end{figure}

Here, we place two actuators at the tip nodes of the wing, as shown in Fig.~\ref{fig:NACAWingSchematic}. We also take the associated two degrees of freedom as observable. With such, we compute DCgains and MHSVs for the first 50 pairs of modes of the wing and plot them in Fig.~\ref{fig:NACAWingHSVDCgain}, from which we see that it is sufficient to take the first pair of modes for the purpose of linear reduction because its normalized DCgain is 96\% and MHSV is 97\%.

\begin{figure}[!ht]
    \centering
    \includegraphics[width=1\linewidth]{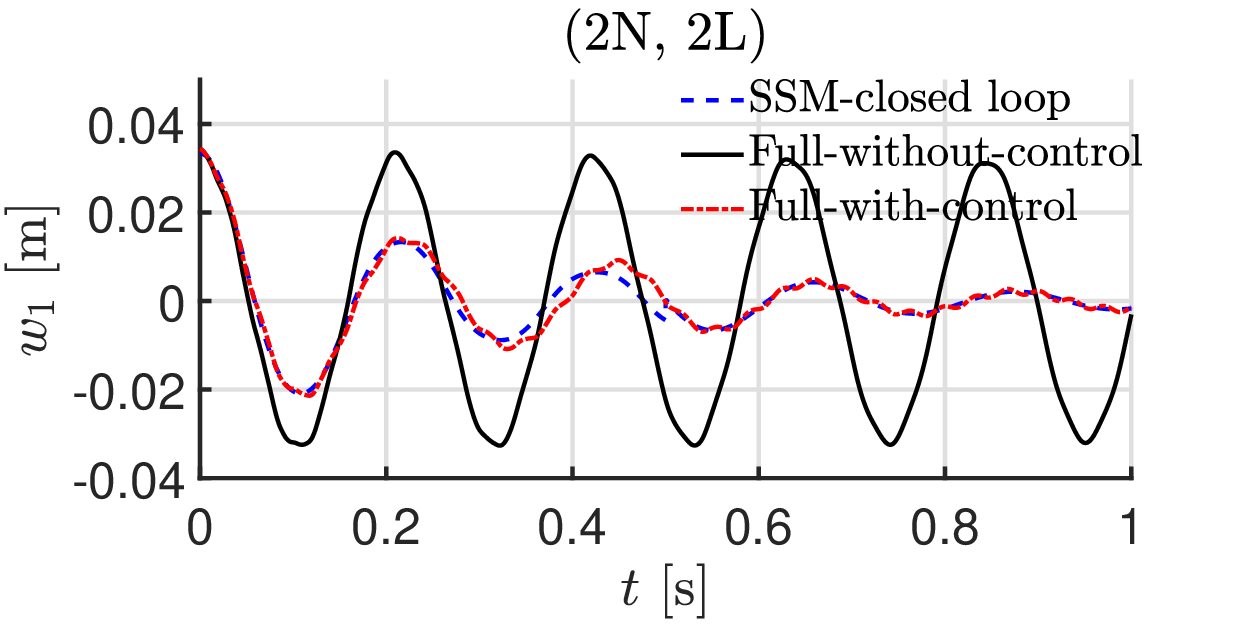}
    \includegraphics[width=1\linewidth]{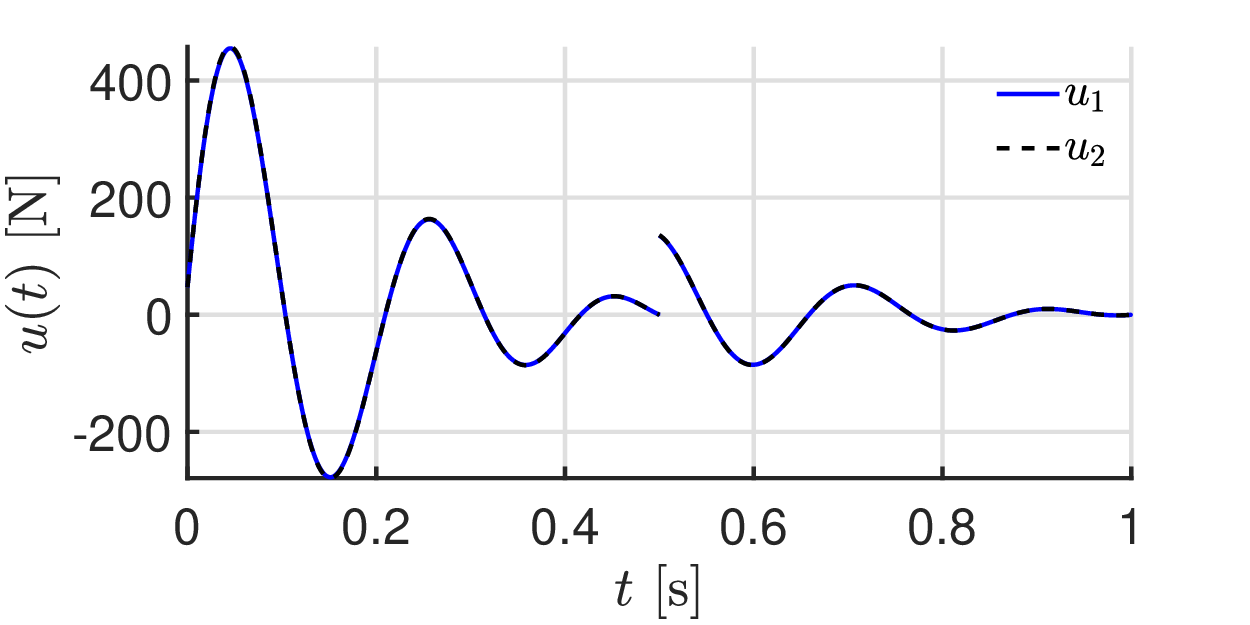}
    \caption{Vibration control of the NACA wing with two actuators placed at two tip nodes of the wing, as shown in Fig.~\ref{fig:NACAWingSchematic}. The time horizon [0,1] is divided into two subintervals with an intersection point at $t=0.5$. In the upper panel, the deflection at tip node 1 is given, where the black solid line denotes the response of the full system without control, the blue dashed line represents the controlled response predicted via SSM-based model reduction, and the red dashed-dotted line gives the response of the full system with the control applied. In the lower panel, the control signals of the two actuators that were obtained using our SSM-based model reduction are given.}
    \label{fig:NACAWing}
\end{figure}

We apply a transverse static loading 2000 N at both tip nodes, resulting in a deflection of 0.0344 $\mathrm{m}$ at the tip nodes. Indeed, this deflection is large enough so that geometric nonlinearity plays an important role~\cite{jain2022compute}. As seen in the upper panel of Fig.~\ref{fig:NACAWing}, the wing undergoes a large amplitude vibration with a slow decay rate. Next, we apply our control design framework to suppress this nonlinear vibration. We set $\epsilon=0.1$, and the control weight matrices are chosen below
\begin{equation}
    \mathbf{Q}=10^4
    \begin{bmatrix}
        \mathbf{I}_n & \mathbf{0} \\
        \mathbf{0}  & \mathbf{I}_n 
    \end{bmatrix},
    \mathbf{\hat{R}}=
    \begin{bmatrix}
        100 & {0} \\
        {0}  & 100
    \end{bmatrix},
    \mathbf{\hat{M}}=\mathbf{0}.
\end{equation}

Moreover, we divide the time horizon $[0,1]$ into two subintervals. As seen in Fig.~\ref{fig:NACAWing}, the large amplitude oscillation is suppressed quickly with our designed control inputs, and more importantly, the controlled response predicted via our SSM-based prediction matches well with that of the full system's response. We observe from the lower panel of Fig.~\ref{fig:NACAWing} that the two control signals are almost the same, which can be explained by the fact that the wing vibrates as a cantilever beam.

\subsection{A fluttered pipe conveying fluid}
\begin{figure}[!ht]
    \centering
    \includegraphics[width=1\linewidth]{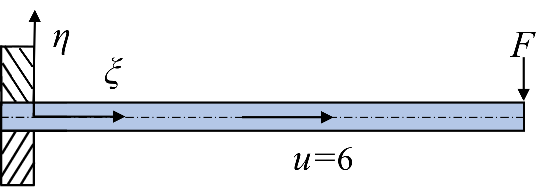}
    \caption{A schematic of a cantilevered pipe conveying fluid}
    \label{fig:FluidPipeSchematic}
\end{figure}

As our last example, we consider the control stabilization of a fluttered pipe conveying fluid shown in Fig~\ref{fig:FluidPipeSchematic}. This example is different from the previous examples in two perspectives. Firstly, we suppress free vibration released from statically deformed configurations in all previous examples. Here, the pipe undergoes a limit cycle oscillation, and we aim to stabilize the flutter motion via active control. Secondly, mechanical systems in previous examples have symmetric stiffness matrix, Rayleigh damping, and only displacement-dependent nonlinear internal forces. The pipe system, however, consists of gyroscopic and follower forces; hence, the associated damping and stiffness matrices are asymmetric~\cite{liNonlinearModelReduction2023a,paidoussis1998fluid}. Moreover, it admits nonlinear damping~\cite{liNonlinearModelReduction2023a,paidoussis1998fluid}.

The governing equation of the pipe can be found in~\cite{liNonlinearModelReduction2023a}. The equation is in the form of a partial differential equation and is discretized into a system of ordinary differential equations (ODEs) via a four-mode Galerkin truncation with eigenfunctions of a cantilever beam, resulting in a system with 4 degrees of freedoms. We refer to~\cite{liNonlinearModelReduction2023a} for more details on the ODEs. In the following computations, we set the dimensionless flow velocity to be 6 such that the pipe is in the post-flutter regime~\cite{liNonlinearModelReduction2023a}.

\begin{figure}[!ht]
    \centering
    \includegraphics[width=1\linewidth]{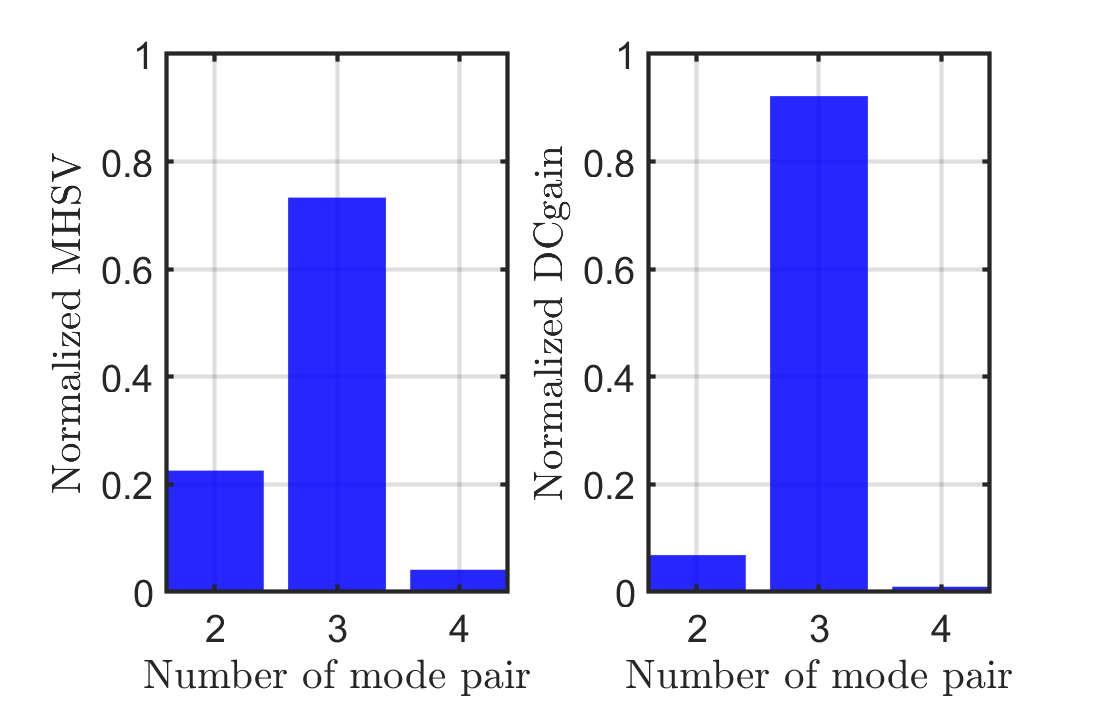}
    \caption{Normalized DCgain and MHSV of the stable eigenmodes of the pipe conveying fluid.}
    \label{fig:FluidDCgain}
\end{figure}

We are concerned with the post-flutter dynamics, and hence, it is sufficient to take the eigenmode associated with eigenvalues $\lambda_{1,2}=0.76\pm13.53\mathrm{i}$ as the master spectral subspace. We compute the associated autonomous unstable SSM at $\mathcal{O}(9)$ expansion. As for the linear reduction, we follow remark~\ref{rk:mhsv}, take the unstable mode into the basis, and compute DCgains and MHSVs for stable modes. Here, we place an actuator at the free end of the pipe and also take the deflection at the free end as observable. The associated DCgains and MHSVs for the three pairs of stable modes are presented in Fig.~\ref{fig:FluidDCgain}. Therefore, we use the first pair of unstable complex modes and the third pair of stable complex modes to reduce the linear part of the system.


\begin{figure}[!ht]
    \centering
    \includegraphics[width=1.05\linewidth]{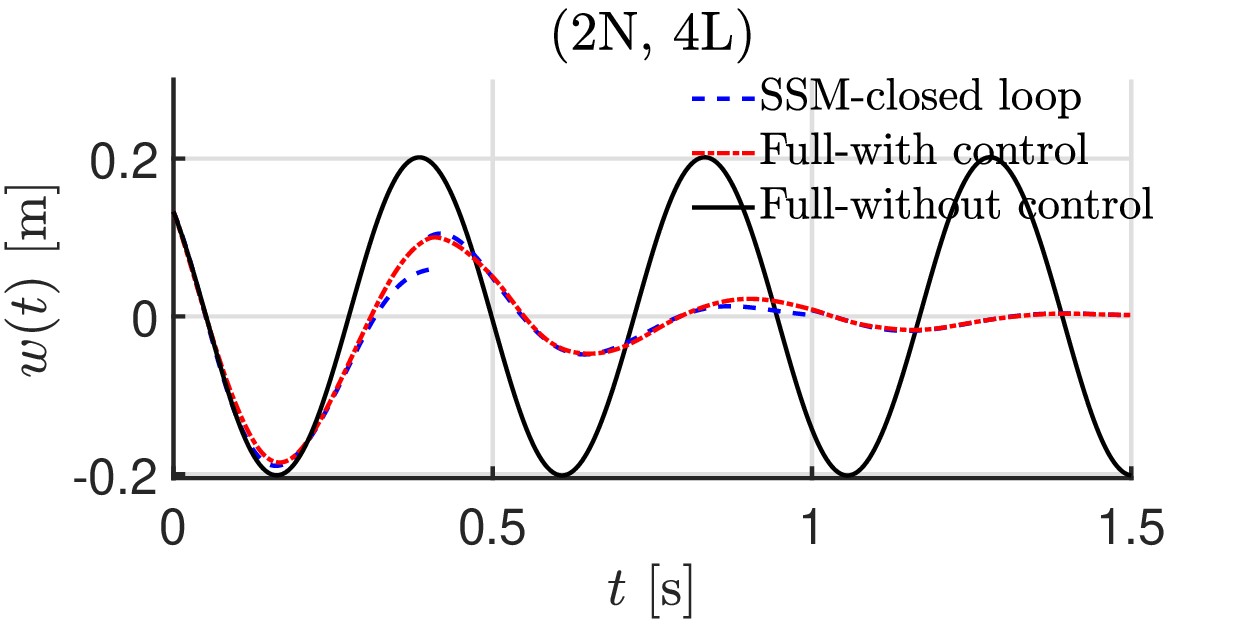}
    \includegraphics[width=1.05\linewidth]{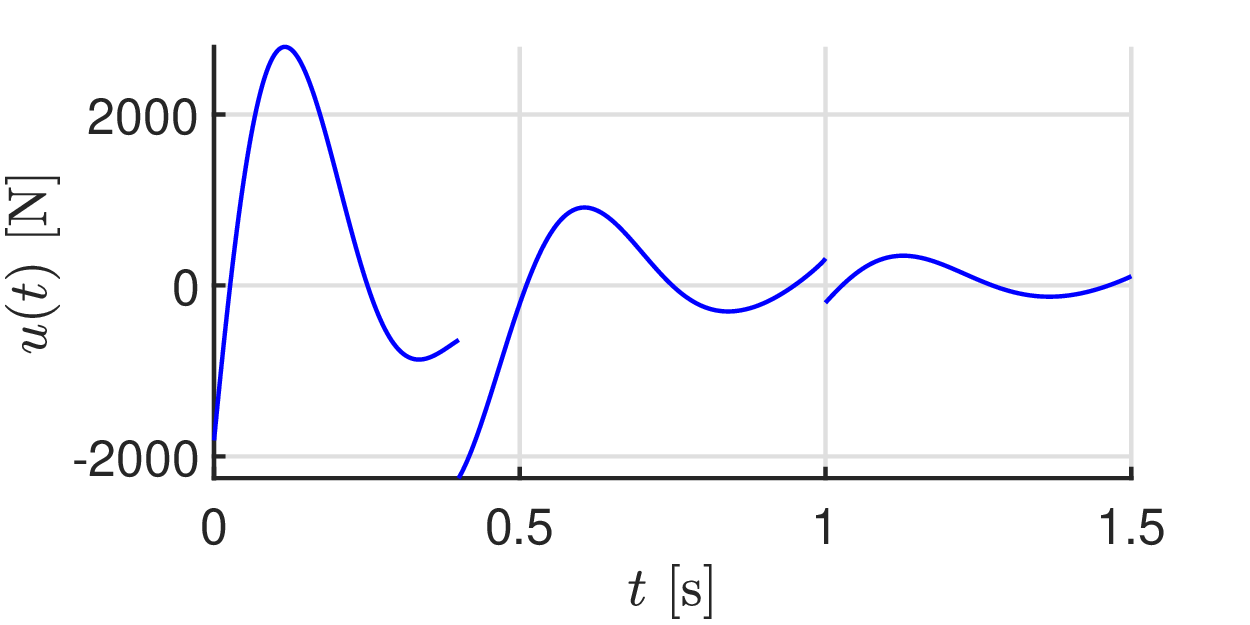}
    \caption{Vibration control of a fluttered pipe conveying fluid with an actuator placed at the free end of the pipe. The time horizon is divided into three subintervals, with two state updates at $t=0.4$ s and $1$ s. In the upper panel, the deflection at the free end is given, where the black solid line denotes the response of the full system without control, the blue dashed line represents the controlled response predicted via SSM-based model reduction, and the red dashed-dotted line gives the response of the full system with the control applied. In the lower panel, the actuator's control signal obtained using our SSM-based model reduction is given.}
    \label{fig:FluidLQR3segs}
\end{figure}

As seen in Fig.~\ref{fig:FluidLQR3segs}, the pipe undergoes a limit cycle oscillation in free vibration. We apply our SSM-based control design framework to suppress the limit cycle motion. Specifically, we divide the time horizon $[0,1.5]$ into three subintervals and choose weight matrices as below
\begin{align}
    & \mathbf{Q}=1.6\times10^5
    \begin{bmatrix}
        \mathbf{I}_n & \mathbf{0} \\
        \mathbf{0}  & \mathbf{I}_n 
    \end{bmatrix},
    \mathbf{\hat{R}}=0.016,  \nonumber \\
    & \mathbf{\hat{M}}=1.6\times10^5
    \begin{bmatrix}
        \mathbf{I}_n & \mathbf{0} \\
        \mathbf{0}  & \mathbf{I}_n 
    \end{bmatrix}.
\end{align}

The designed control signals for the three subintervals are shown in the lower panel of Fig.~\ref{fig:FluidLQR3segs}. We observe from the upper panel of the figure that the system is stabilized to equilibrium after about 3 cycles of vibration, and the SSM-based prediction matches well with that of the response of the full system. This again demonstrates the effectiveness of our SSM-based control design framework.

\section{Conclusions}
\label{sec:conclusion}

In this study, we use spectral submanifolds (SSMs) to construct low-dimensional reduced-order models (ROMs) that facilitate the real-time control of large amplitude vibration of flexible structures. We transform high-dimensional nonlinear optimal control problems into low-dimensional linear optimal control problems, which enables the use of well-established linear control theory and effectively solves the reduced control problems. We also achieve a significant dimension reduction because our SSM-based ROMs are often two-dimensional for systems without internal resonances, paving the road for the real-time control of the nonlinear vibrations of flexible structures.

The established control design framework via SSM-based model reduction displays wide applicability, including systems possibly with external forcing, internal resonance, gyroscopic and follower forces, and nonlinear damping. Our SSM-based control design can suppress the nonlinear vibrations of these systems effectively. In particular, we have successfully suppressed the nonlinear vibration of an aircraft wing modeled using the finite element method with more than 130,000 degrees of freedom. We have also effectively stabilized the flutter motion of a cantilevered pipe conveying fluid.

In short, our control design framework via SSM-based model reduction provides a powerful approach for the active control of nonlinear vibrations of flexible structures modeled via high-dimensional discrete models. The demonstrated advantages and benefits of the SSM-based reduction can be carried over to other tasks, such as trajectory tracking of soft robotics. Our SSM-based control design framework will play an important role in the vibration control design via digital twins for more complicated engineering applications.

\backmatter

\bmhead{Data availability}
Data will be made available on request.

\bmhead{Acknowledgements}

The financial support of the National Natural Science Foundation of China (No. 12302014) and Shenzhen Science and Technology Innovation Commission (No. 20231115172355001) for this work is gratefully acknowledged.


\section*{Conflict of interest}
The authors declare that they have no known competing financial interests or personal relationships that could have appeared to influence the work reported in this paper.

\makeatletter
\makeatother

\begin{thebibliography}{72}
\providecommand{\natexlab}[1]{#1}
\providecommand{\url}[1]{{#1}}
\providecommand{\urlprefix}{URL }
\providecommand{\doi}[1]{\url{https://doi.org/#1}}
\providecommand{\eprint}[2][]{\url{#2}}
 \bibcommenthead

\bibitem[{Wani et~al(2022)Wani, Tantray, Noroozinejad~Farsangi, Nikitas, Noori, Samali, and Yang}]{waniCriticalReviewControl2022}
Wani ZR, Tantray M, Noroozinejad~Farsangi E, et~al (2022) A {{Critical Review}} on {{Control Strategies}} for {{Structural Vibration Control}}. Annual Reviews in Control 54:103--124

\bibitem[{Wagg and Neild(2010)}]{davidwagg}
Wagg D, Neild S (2010) Nonlinear Vibration with Control for Flexible and Adaptive Structures. Springer

\bibitem[{Xie and Aly(2020)}]{xieStructuralControlVibration2020}
Xie F, Aly AM (2020) Structural control and vibration issues in wind turbines: {{A}} review. Engineering Structures 210:110087

\bibitem[{Webster et~al(2009)Webster, Romano, and Cowan}]{websterMechanicsPrecurvedTubeContinuum2009}
Webster RJ, Romano JM, Cowan NJ (2009) Mechanics of {{Precurved-Tube Continuum Robots}}. IEEE Transactions on Robotics 25(1):67--78

\bibitem[{Russo et~al(2023)Russo, Sadati, Dong, Mohammad, Walker, Bergeles, Xu, and Axinte}]{russoContinuumRobotsOverview2023}
Russo M, Sadati SMH, Dong X, et~al (2023) Continuum {{Robots}}: {{An Overview}}. Advanced Intelligent Systems p 2200367

\bibitem[{{Lee-Glauser} et~al(1996){Lee-Glauser}, Ahmadi, and Layton}]{lee-glauserSatelliteActivePassive1996}
{Lee-Glauser} GJ, Ahmadi G, Layton JB (1996) Satellite active and passive vibration control during liftoff. Journal of Spacecraft and Rockets 33(3):428--432

\bibitem[{Mohamed et~al(2005)Mohamed, Martins, Tokhi, S{\'a}~Da~Costa, and Botto}]{mohamedVibrationControlVery2005}
Mohamed Z, Martins J, Tokhi M, et~al (2005) Vibration control of a very flexible manipulator system. Control Engineering Practice 13(3):267--277

\bibitem[{{El-Khoury} and Adeli(2013)}]{khouryRecentAdvancesVibration2013}
{El-Khoury} O, Adeli H (2013) Recent {{Advances}} on {{Vibration Control}} of {{Structures Under Dynamic Loading}}. Archives of Computational Methods in Engineering 20(4):353--360

\bibitem[{Ghaedi et~al(2017)Ghaedi, Ibrahim, Adeli, and Javanmardi}]{ghaediInvitedReviewRecent2017}
Ghaedi K, Ibrahim Z, Adeli H, et~al (2017) Invited {{Review}}: {{Recent}} developments in vibration control of building and bridge structures. Journal of Vibroengineering 19(5):3564--3580

\bibitem[{Kumar et~al(2023)Kumar, Kumar, and Kumar}]{kumarReviewControllersStructural2023}
Kumar G, Kumar R, Kumar A (2023) A {{Review}} of the {{Controllers}} for {{Structural Control}}. Archives of Computational Methods in Engineering 30(6):3977--4000

\bibitem[{Li and Zhu(2022)}]{liSelfPoweredActiveVibration2022}
Li JY, Zhu S (2022) Self-{{Powered Active Vibration Control}}: {{Concept}}, {{Modeling}}, and {{Testing}}. Engineering 11:126--137

\bibitem[{Lu et~al(2004)Lu, Chung, and Lin}]{luGeneralMethodSemiActive2004}
Lu LY, Chung LL, Lin GL (2004) A {{General Method}} for {{Semi-Active Feedback Control}} of {{Variable Friction Dampers}}. Journal of Intelligent Material Systems and Structures 15(5):393--412

\bibitem[{Shi et~al(2017)Shi, Zhu, and Nagarajaiah}]{shiPerformanceComparisonPassive2017}
Shi X, Zhu S, Nagarajaiah S (2017) Performance {{Comparison}} between {{Passive Negative-Stiffness Dampers}} and {{Active Control}} in {{Cable Vibration Mitigation}}. Journal of Bridge Engineering 22(9):04017054

\bibitem[{Thenozhi and Yu(2013)}]{thenozhiAdvancesModelingVibration2013}
Thenozhi S, Yu W (2013) Advances in modeling and vibration control of building structures. Annual Reviews in Control 37(2):346--364

\bibitem[{Korkmaz(2011)}]{korkmazReviewActiveStructural2011}
Korkmaz S (2011) A review of active structural control: Challenges for engineering informatics. Computers \& Structures 89(23-24):2113--2132

\bibitem[{Zhang et~al(2017)Zhang, Zang, Li, Wang, and Li}]{zhangActivepassiveIntegratedVibration2017}
Zhang Y, Zang Y, Li M, et~al (2017) Active-passive integrated vibration control for control moment gyros and its application to satellites. Journal of Sound and Vibration 394:1--14

\bibitem[{Ding and Ji(2023)}]{dingVibrationControlFluidconveying2023}
Ding H, Ji JC (2023) Vibration control of fluid-conveying pipes: A state-of-the-art review. Applied Mathematics and Mechanics 44(9):1423--1456

\bibitem[{Thenozhi and Yu(2014)}]{thenozhiStabilityAnalysisActive2014}
Thenozhi S, Yu W (2014) Stability analysis of active vibration control of building structures using {{PD}}/{{PID}} control. Engineering Structures 81:208--218

\bibitem[{Ho and Ma(2007)}]{hoActiveVibrationControl2007}
Ho CC, Ma CK (2007) Active vibration control of structural systems by a combination of the linear quadratic {{Gaussian}} and input estimation approaches. Journal of Sound and Vibration 301(3-5):429--449

\bibitem[{Takamoto et~al(2022)Takamoto, Abe, Hara, Otsuka, and Makihara}]{takamotoComprehensivePredictiveControl2022}
Takamoto I, Abe M, Hara Y, et~al (2022) Comprehensive predictive control for vibration suppression based on piecewise constant input formulation. Journal of Intelligent Material Systems and Structures 33(7):901--917

\bibitem[{Canciello and Cavallo(2017)}]{cancielloSelectiveModalControl2017}
Canciello G, Cavallo A (2017) Selective modal control for vibration reduction in flexible structures. Automatica 75:282--287

\bibitem[{Choi and Han(2007)}]{choiVibrationControlElectrorheological2007}
Choi SB, Han YM (2007) Vibration control of electrorheological seat suspension with human-body model using sliding mode control. Journal of Sound and Vibration 303(1-2):391--404

\bibitem[{Li and Adeli(2018)}]{liControlMethodologiesVibration2018}
Li Z, Adeli H (2018) Control methodologies for vibration control of smart civil and mechanical structures. Expert Systems 35(6):e12354

\bibitem[{Kandasamy et~al(2016)Kandasamy, Cui, Townsend, Foo, Guo, Shenoi, and Xiong}]{kandasamyReviewVibrationControl2016}
Kandasamy R, Cui F, Townsend N, et~al (2016) A review of vibration control methods for marine offshore structures. Ocean Engineering 127:279--297

\bibitem[{Gao et~al(2021)Gao, Yu, Zhang, Wang, and Zhai}]{gaoVibrationAnalysisControl2021}
Gao P, Yu T, Zhang Y, et~al (2021) Vibration analysis and control technologies of hydraulic pipeline system in aircraft: {{A}} review. Chinese Journal of Aeronautics 34(4):83--114

\bibitem[{Jian(2015)}]{SHLX201504001}
Jian X (2015) Advances of research on vibration control. Chinese quarterly of mechanics 36(04):547--565

\bibitem[{Kangunde et~al(2021)Kangunde, Jamisola, and Theophilus}]{kangundeReviewDronesControlled2021}
Kangunde V, Jamisola RS, Theophilus EK (2021) A review on drones controlled in real-time. International Journal of Dynamics and Control 9(4):1832--1846

\bibitem[{Jain and Haller(2022)}]{jain2022compute}
Jain S, Haller G (2022) How to compute invariant manifolds and their reduced dynamics in high-dimensional finite element models. Nonlinear dynamics 107(2):1417--1450

\bibitem[{Wan et~al(2024)Wan, Ma, Dong, Luo, and Ni}]{wanDatadrivenModelReduction2024}
Wan HP, Ma Q, Dong GS, et~al (2024) Data-driven model reduction approach for active vibration control of cable-strut structures. Engineering Structures 302:117434

\bibitem[{Besselink et~al(2013)Besselink, Tabak, Lutowska, Van De~Wouw, Nijmeijer, Rixen, Hochstenbach, and Schilders}]{besselinkComparisonModelReduction2013}
Besselink B, Tabak U, Lutowska A, et~al (2013) A comparison of model reduction techniques from structural dynamics, numerical mathematics and systems and control. Journal of Sound and Vibration 332(19):4403--4422

\bibitem[{Suman and Kumar(2022)}]{sumanInvestigationImplementationModel2022}
Suman SK, Kumar A (2022) Investigation and {{Implementation}} of {{Model Order Reduction Technique}} for {{Large Scale Dynamical Systems}}. Archives of Computational Methods in Engineering 29(5):3087--3108

\bibitem[{Xianmin et~al(2002)Xianmin, Changjian, and Erdman}]{xianminActiveVibrationController2002}
Xianmin Z, Changjian S, Erdman AG (2002) Active vibration controller design and comparison study of flexible linkage mechanism systems. Mechanism and Machine Theory 37(9):985--997

\bibitem[{Vakilzadeh et~al(2020)Vakilzadeh, Vatankhah, and Eghtesad}]{vakilzadehVibrationControlMicroscale2020}
Vakilzadeh M, Vatankhah R, Eghtesad M (2020) Vibration control of micro-scale structures using their reduced second order bilinear models based on multi-moment matching criteria. Applied Mathematical Modelling 78:287--296

\bibitem[{Gildin et~al(2009)Gildin, Antoulas, Sorensen, and Bishop}]{gildinModelControllerReduction2009a}
Gildin E, Antoulas AC, Sorensen D, et~al (2009) Model and controller reduction applied to structural control using passivity theory. Structural Control and Health Monitoring 16(3):319--334

\bibitem[{Banks et~al(2002)Banks, Del~Rosario, and Tran}]{banksProperOrthogonalDecompositionbased2002}
Banks H, Del~Rosario R, Tran H (2002) Proper orthogonal decomposition-based control of transverse beam vibrations: Experimental implementation. IEEE Transactions on Control Systems Technology 10(5):717--726

\bibitem[{Mathews et~al(2002)Mathews, Sule, and Venkatesan}]{mathewsOrderReductionClosedLoop2002}
Mathews A, Sule VR, Venkatesan C (2002) Order {{Reduction}} and {{Closed-Loop Vibration Control}} in {{Helicopter Fuselages}}. Journal of Guidance, Control, and Dynamics 25(2):316--323

\bibitem[{Antoulas(2005)}]{antoulasApproximationLargescaleDynamical2005}
Antoulas AC (2005) Approximation of Large-Scale Dynamical Systems. Advances in Design and Control, {Society for Industrial and Applied Mathematics}, Philadelphia

\bibitem[{King et~al(2006)King, Hovakimyan, Evans, and Buhl}]{kingReducedOrderControllers2006}
King BB, Hovakimyan N, Evans KA, et~al (2006) Reduced order controllers for distributed parameter systems: {{LQG}} balanced truncation and an adaptive approach. Mathematical and Computer Modelling 43(9-10):1136--1149

\bibitem[{Zhou and Chen(2015)}]{zhouIntelligentVibrationControl2015}
Zhou L, Chen G (2015) Intelligent {{Vibration Control}} for {{High-Speed Spinning Beam Based}} on {{Fuzzy Self-Tuning PID Controller}}. Shock and Vibration 2015:1--8

\bibitem[{Haller and Ponsioen(2016)}]{hallerNonlinearNormalModes2016a}
Haller G, Ponsioen S (2016) Nonlinear normal modes and spectral submanifolds: Existence, uniqueness and use in model reduction. Nonlinear Dynamics 86(3):1493--1534

\bibitem[{Touz{\'e} et~al(2021)Touz{\'e}, Vizzaccaro, and Thomas}]{touzeModelOrderReduction2021}
Touz{\'e} C, Vizzaccaro A, Thomas O (2021) Model order reduction methods for geometrically nonlinear structures: A review of nonlinear techniques. Nonlinear Dynamics 105(2):1141--1190

\bibitem[{Ma et~al(2024)Ma, Wang, Chen, Brighton, and Anish}]{maSuspensionNonlinearAnalysis2024}
Ma G, Wang P, Chen L, et~al (2024) Suspension nonlinear analysis and {{VSS-LMS}} adaptive filtering control of satellite borne flexible structure. Nonlinear Dynamics 112(5):3679--3693

\bibitem[{Nechak(2022)}]{nechakRobustNonlinearControl2022}
Nechak L (2022) Robust nonlinear control synthesis by using centre manifold-based reduced models for the mitigating of friction-induced vibration. Nonlinear Dynamics 108(3):1885--1901

\bibitem[{Shaw and Pierre(1999)}]{shaw1999modal}
Shaw SW, Pierre C (1999) Modal analysis-based reduced-order models for nonlinear structures andmdash; an invariant manifold approach. The shock and vibration digest 31(1):3--16

\bibitem[{Touz{\'e} and Amabili(2006)}]{touze2006nonlinear}
Touz{\'e} C, Amabili M (2006) Nonlinear normal modes for damped geometrically nonlinear systems: Application to reduced-order modelling of harmonically forced structures. Journal of Sound and Vibration 298(4-5):958--981

\bibitem[{Liu and Wagg(2019)}]{Liu2019SimultaneousNF}
Liu X, Wagg DJ (2019) Simultaneous normal form transformation and model-order reduction for systems of coupled nonlinear oscillators. Proceedings of the Royal Society A 475

\bibitem[{Li and Haller(2022)}]{liNonlinearAnalysisForced2022}
Li M, Haller G (2022) Nonlinear analysis of forced mechanical systems with internal resonance using spectral submanifolds, {{Part II}}: {{Bifurcation}} and quasi-periodic response. Nonlinear Dynamics 110(2):1045--1080

\bibitem[{Li et~al(2023)Li, Jain, and Haller}]{liModelReductionConstrained2023}
Li M, Jain S, Haller G (2023) Model reduction for constrained mechanical systems via spectral submanifolds. Nonlinear Dynamics 111(10):8881--8911

\bibitem[{Cenedese et~al(2022)Cenedese, Ax{\aa}s, B{\"a}uerlein, Avila, and Haller}]{cenedese2022data}
Cenedese M, Ax{\aa}s J, B{\"a}uerlein B, et~al (2022) Data-driven modeling and prediction of non-linearizable dynamics via spectral submanifolds. Nature Communications 13(1):1--13

\bibitem[{Kasz{\'a}s and Haller(2024)}]{kaszas2024capturing}
Kasz{\'a}s B, Haller G (2024) Capturing the edge of chaos as a spectral submanifold in pipe flows. Journal of Fluid Mechanics 979:A48

\bibitem[{Li et~al(2023)Li, Yan, and Wang}]{liNonlinearModelReduction2023a}
Li M, Yan H, Wang L (2023) Nonlinear model reduction for a cantilevered pipe conveying fluid: {{A}} system with asymmetric damping and stiffness matrices. Mechanical Systems and Signal Processing 188:109993

\bibitem[{Li et~al(2024)Li, Yan, and Wang}]{li2024data}
Li M, Yan H, Wang L (2024) Data-driven model reduction for pipes conveying fluid via spectral submanifolds. International Journal of Mechanical Sciences p 109414

\bibitem[{Xu et~al(2024)Xu, Kaszás, Cenedese, Berti, Coletti, and Haller}]{Xu_Kaszás_Cenedese_Berti_Coletti_Haller_2024}
Xu Z, Kaszás B, Cenedese M, et~al (2024) Data-driven modelling of the regular and chaotic dynamics of an inverted flag from experiments. Journal of Fluid Mechanics 987:R7

\bibitem[{Breunung and Haller(2018)}]{breunung2018explicit}
Breunung T, Haller G (2018) Explicit backbone curves from spectral submanifolds of forced-damped nonlinear mechanical systems. Proceedings of the Royal Society A: Mathematical, Physical and Engineering Sciences 474(2213):20180083

\bibitem[{Li et~al(2024)Li, Jain, and Haller}]{li2024fast}
Li M, Jain S, Haller G (2024) Fast computation and characterization of forced response surfaces via spectral submanifolds and parameter continuation. Nonlinear Dynamics pp 1--27

\bibitem[{Haller and Kaundinya(2024)}]{Haller2024NonlinearMR}
Haller G, Kaundinya RS (2024) Nonlinear model reduction to temporally aperiodic spectral submanifolds. Chaos 34:043152

\bibitem[{Alora et~al(2023{\natexlab{a}})Alora, Cenedese, Schmerling, Haller, and Pavone}]{alora2023data}
Alora JI, Cenedese M, Schmerling E, et~al (2023{\natexlab{a}}) Data-driven spectral submanifold reduction for nonlinear optimal control of high-dimensional robots. In: 2023 IEEE International Conference on Robotics and Automation (ICRA), IEEE, pp 2627--2633

\bibitem[{Alora et~al(2023{\natexlab{b}})Alora, Pabon, K{\"o}hler, Cenedese, Schmerling, Zeilinger, Haller, and Pavone}]{alora2023robust}
Alora JI, Pabon LA, K{\"o}hler J, et~al (2023{\natexlab{b}}) Robust nonlinear reduced-order model predictive control. In: 2023 62nd IEEE Conference on Decision and Control (CDC), IEEE, pp 4798--4805

\bibitem[{Haller and Ponsioen(2016)}]{haller2016nonlinear}
Haller G, Ponsioen S (2016) Nonlinear normal modes and spectral submanifolds: existence, uniqueness and use in model reduction. Nonlinear Dynamics 86(3):1493--1534

\bibitem[{Haller et~al(2023)Haller, Kasz{\'a}s, Liu, and Ax{\aa}s}]{haller2023nonlinear}
Haller G, Kasz{\'a}s B, Liu A, et~al (2023) Nonlinear model reduction to fractional and mixed-mode spectral submanifolds. Chaos: An Interdisciplinary Journal of Nonlinear Science 33(6)

\bibitem[{Li et~al(2022)Li, Jain, and Haller}]{part-i}
Li M, Jain S, Haller G (2022) Nonlinear analysis of forced mechanical systems with internal resonance using spectral submanifolds, {Part I}: Periodic response and forced response curve. Nonlinear Dynamics 110:1005--1043

\bibitem[{Jain et~al(2023)Jain, Thurnher, Li, and Haller}]{ssmtool21}
Jain S, Thurnher T, Li M, et~al (2023) {SSMTool} 2.5: Computation of invariant manifolds \& their reduced dynamics in high-dimensional mechanics problems. \url{https://doi.org/10.5281/zenodo.10018285}, accessed: 2024-7-1

\bibitem[{Rantzer(2011)}]{rantzerDistributedControlPositive2011}
Rantzer A (2011) Distributed control of positive systems. In: {{IEEE Conference}} on {{Decision}} and {{Control}} and {{European Control Conference}}. IEEE, Orlando, FL, USA, pp 6608--6611

\bibitem[{Kawano et~al(2020)Kawano, Besselink, Scherpen, and Cao}]{kawanoDataDrivenModelReduction2020}
Kawano Y, Besselink B, Scherpen JMA, et~al (2020) Data-{{Driven Model Reduction}} of {{Monotone Systems}} by {{Nonlinear DC Gains}}. IEEE Transactions on Automatic Control 65(5):2094--2106

\bibitem[{Chang et~al(2002)Chang, Gopinathan, Varadan, and Varadan}]{Chang2002DesignOR}
Chang W, Gopinathan SV, Varadan VV, et~al (2002) Design of robust vibration controller for a smart panel using finite element model. Journal of Vibration and Acoustics 124:265--276

\bibitem[{Chang(2004)}]{changModelReductionBased2004}
Chang W (2004) Model reduction based on modal {{Hankel}} singular values. In: Smith RC (ed) Smart {{Structures}} and {{Materials}}, San Diego, CA, p 433

\bibitem[{Kung and Lin(1980)}]{Kung1980OptimalHM}
Kung S, Lin DW (1980) Optimal hankel-norm model reductions: Multivariable systems. 1980 19th IEEE Conference on Decision and Control including the Symposium on Adaptive Processes pp 187--194

\bibitem[{Moore(1981)}]{Moore1981PrincipalCA}
Moore BA (1981) Principal component analysis in linear systems: Controllability, observability, and model reduction. IEEE Transactions on Automatic Control 26:17--32

\bibitem[{Mattingley et~al(2011)Mattingley, Wang, and Boyd}]{Mattingley2011RecedingHC}
Mattingley J, Wang Y, Boyd S (2011) Receding horizon control. IEEE Control Systems 31:52--65

\bibitem[{Nayfeh et~al(1974)Nayfeh, Mook, and Sridhar}]{Nayfeh1974NonlinearAO}
Nayfeh AH, Mook DT, Sridhar S (1974) Nonlinear analysis of the forced response of structural elements. Journal of the Acoustical Society of America 55:281--291

\bibitem[{Jain et~al(2017)Jain, Tiso, Rutzmoser, and Rixen}]{jainQuadraticManifoldModel2017}
Jain S, Tiso P, Rutzmoser JB, et~al (2017) A quadratic manifold for model order reduction of nonlinear structural dynamics. Computers \& Structures 188:80--94

\bibitem[{Pa{\"\i}doussis(1998)}]{paidoussis1998fluid}
Pa{\"\i}doussis MP (1998) Fluid-structure interactions: slender structures and axial flow, vol~1. Academic press

\end{thebibliography}


\end{document}